\documentclass[12pt]{article}
\usepackage{makeidx}
\usepackage{amssymb,amsthm,amsmath,amsfonts,enumerate}
\usepackage{setspace,graphicx,caption,subcaption,float,relsize}
\usepackage{url}
\usepackage{hyperref}
\usepackage{rotating}
\usepackage[authoryear,longnamesfirst]{natbib}
\usepackage{array}
\usepackage{indentfirst}

\hypersetup{
  colorlinks = true, urlcolor = cyan, linkcolor = blue, citecolor = red }
\renewenvironment{abstract}
 {\small
  \begin{center}
  \bfseries \abstractname\vspace{-.0em}\vspace{0pt}
  \end{center}
  \list{}{    \setlength{\leftmargin}{0mm}
    \setlength{\rightmargin}{\leftmargin}  }  \item\relax}
 {\endlist}
\makeatletter
\def\maketag@@@#1{\hbox{\m@th\normalfont\normalsize#1}}
\makeatother

\allowdisplaybreaks
\newtheorem {theorem}{Theorem}[section]
\newtheorem {assumption}{Assumption}

\newtheorem{lemma}[theorem]{Lemma}

\newtheorem{proposition}[theorem]{Proposition}
\newtheorem{remark}{Remark}[section]

\allowdisplaybreaks

\counterwithin*{equation}{section}
\numberwithin{equation}{section}

\begin{document}

\title{$\mathcal{L}_{q}$-maximal inequality for high dimensional \\
means under dependence}
\author{Jonathan B. Hill\thanks{%
Department of Economics, University of North Carolina, Chapel Hill, North
Carolina, E-mail:\texttt{jbhill@email.unc.edu}; \texttt{%
https://tarheels.live/jbhill}.}\medskip \\
Dept. of Economics, University of North Carolina, Chapel Hill, NC}
\date{{\large This draft:} \today
}
\maketitle

\begin{abstract}
We derive an $\mathcal{L}_{q}$-maximal inequality for zero mean dependent
random variables $\{x_{t}\}_{t=1}^{n}$ on $\mathbb{R}^{p}$, where $p$ $>>$ $%
n $ is allowed. The upper bound is a familiar multiple of $\ln (p)$ and an $%
l_{\infty }$ moment, as well as Kolmogorov distances based on Gaussian
approximations $(\rho _{n},\tilde{\rho}_{n})$, derived with and without
negligible truncation and sub-sample blocking. The latter arise due to a
departure from independence and therefore a departure from standard
symmetrization arguments. Examples are provided demonstrating $(\rho _{n},%
\tilde{\rho}_{n})$ $\rightarrow $ $0$ under heterogeneous mixing and
physical dependence conditions, where $(\rho _{n},\tilde{\rho}_{n})$ are
multiples of $\ln (p)/n^{b}$ for some $b$ $>$ $0$ that depends on memory,
tail decay, the truncation level and block size.\medskip \newline
 \smallskip \newline
\end{abstract}

\setstretch{1.2}

\section{Introduction\label{sec:intro}}

In this paper we establish an $\mathcal{L}_{q}$-maximal inequality for zero
mean, high dimensional dependent random variables $x_{t}$ $=$ $%
[x_{i,t}]_{i=1}^{p}$ on $\mathbb{R}^{p}$, where $p$ $>>$ $n$ and $p$ $=$ $%
p_{n}$ $\rightarrow $ $\infty $ as $n$ $\rightarrow $ $\infty $ are allowed.
Let $\{x_{t}\}_{t=1}^{n}$ be a sample of observations, $n$ $\in $ $\mathbb{N}
$, and assume $\max_{1\leq i\leq p}\max_{1\leq t\leq n}|x_{i,t}|$ is $%
\mathcal{L}_{q}$ integrable for some $q$ $\geq $ $1$. We show $\mathbb{E}%
\max_{1\leq i\leq p}|1/n\sum_{t=1}^{n}x_{i,t}|^{q}$ is bounded by a multiple
of $(\ln \left( p\right) /\mathcal{N}_{n})^{q/2}$ in a general dependence
setting that uses blocking, where $\mathcal{N}_{n}$ $\rightarrow $ $\infty $%
\ are the number of blocks satisfying $\mathcal{N}_{n}=o(n)$. We then prove
required Gaussian approximations and comparisons are asymptotically
negligible under mixing and physical dependence properties. The maximal
moment becomes $o(1)$ after bounding the growth of $p$, depending on whether 
$x_{i,t}$ have sub-exponential tails. We use a negligible truncation
approximation with a diverging level $\mathcal{M}_{n}$, and a blocking (or $%
M $-dependent) approximation popular in the dependent multiplier bootstrap
literature (e.g. \citet{Liu1988,Hansen1996,Shao2010,Shao2011}). Thus an
upper bound on $p$ will naturally depend on truncation $\mathcal{M}_{n}$ and
block size $b_{n}$, as well as dependence decay and heterogeneity. See
below, and Section \ref{sec:Lq_ineq}, for notation conventions.

Under high dimensionality concentration inequalities can be combined with a $%
log$-$exp$ or \textquotedblleft \textit{smooth max}\textquotedblright\
approximation to bound $\mathbb{E}\max_{1\leq i\leq
p}|1/n\sum_{t=1}^{n}x_{i,t}|$. Thus, for example, for any $\lambda $ $>$ $0$
we have by Jensen's inequality 
\begin{equation}
\mathbb{E}\max_{1\leq i\leq p}\left\vert \frac{1}{n}\sum_{t=1}^{n}x_{i,t}%
\right\vert \leq \frac{1}{\lambda }\ln \left( p\max_{1\leq i\leq
p}\int_{0}^{\infty }\mathbb{P}\left( \left\vert \frac{1}{n}%
\sum_{t=1}^{n}x_{i,t}\right\vert >\frac{1}{\lambda }\ln \left( u\right)
\right) du\right) .  \label{log_exp}
\end{equation}%
The integrated probability can be assessed under sub-exponential tails using
probability concentration or large deviation bounds like Bernstein and
Fuk-Nagaev inequalities. The latter have been extended to weakly dependent
data in many mixing and non-mixing contexts 
\citep[see,
e.g.,][]{Rio2017,Hill2024_maxlln}. If tail decay is slower than exponential
then use Lyapunov's inequality and $\max_{1\leq i\leq p}|x_{i}|$ $\leq $ $%
\sum_{i=1}^{p}|x_{i}|$ to deduce $\mathbb{E}{\small {\max_{1\leq i\leq p}|1/n%
}}\sum_{t=1}^{n}x_{i,t}|^{q}$ $\leq $ $(p^{q/s}/n^{q/2})\max_{1\leq i\leq
p}||1/\sqrt{{\small {n}}}\sum\nolimits_{t=1}^{n}x_{i,t}||_{s}^{q}$ for $s$ $%
\geq $ $q$. The final $\mathcal{L}_{q}$-moment is boundable in many
dependence and heterogeneity settings %
\citep[e.g.][]{Davydov1968,McLeish1975,Wu2005}, yielding $\mathbb{E}{\small {%
\max_{1\leq i\leq p}|1/n}}\sum_{t=1}^{n}x_{i,t}|^{q}$ $\rightarrow $ $0$
under a polynomial bound on $p$.

As far as we know, \cite{Nemirovski2000}-type bounds for $\mathbb{E}%
\max_{1\leq i\leq p}|1/n\sum_{t=1}^{n}x_{i,t}|^{q}$ that are a multiple of $%
\ln (p)$ \textit{without} exploiting sub-exponential tails only exist for
independent data. The latter hinges on a symmetrization argument classically
laid out in \citet[Chapt. II.2]{Pollard1984} and 
\citet[Chapt.'s
2.2-2.3]{vanderVaartWellner1996}. See also Section \ref{sec:Lq_ineq} below.
Write $\max_{i}$ $=$ $\max_{1\leq i\leq p}$ and $\max_{i,t}$ $=$ $%
\max_{1\leq i\leq p}\max_{1\leq t\leq n}$, and define $\bar{x}_{i,n}$ $:=$ $%
1/n\sum_{t=1}^{n}x_{i,t}$. In this setting \cite{Nemirovski2000} shows for
any $q$ $\geq $ $1$ and $p$ $\geq $ $e^{q-1} $\ 
\citep[see, e.g.,][Theorem
14.24]{BuhlmannVanDeGeer2011} 
\begin{equation}
\mathbb{E}\max_{i}\left\vert \bar{x}_{i,n}\right\vert ^{q}\leq C_{n}(p,q)%
\mathbb{E}\left( \max_{i}\frac{1}{n}\sum_{t=1}^{n}x_{i,t}^{2}\right) ^{q/2}%
\text{ where }C_{n}(p,q)=\left( \frac{8\ln \left( 2p\right) }{n}\right)
^{q/2}.  \label{Nem_q}
\end{equation}%
The bound depends on $p$ both from $\ln (p)$ and $\mathbb{E}%
(\max_{i}1/n\sum_{t=1}^{n}x_{i,t}^{2})^{q/2}$. If $q$ $=$ $2$ then of course 
$\mathbb{E}\max_{i}\bar{x}_{i,n}^{2}$ $\lesssim $ $\sqrt{\ln (p)/n}%
\sum_{t=1}^{n}\mathbb{E}\max_{i}x_{i,t}^{2}/n$. If $x_{i,t}$ is sub-Gaussian
then $\mathbb{E}\max_{i}x_{i,t}^{2}$ $=$ $O(\ln (n))$; otherwise $\mathbb{E}%
\max_{i}x_{i,t}^{2}$ $\leq $ $p^{2/s}(\max_{i}||x_{i,t}||_{s})^{2}$ under $%
\mathcal{L}_{s}$-boundedness, $s$ $\geq $ $2$, by Lyapunov's inequality and $%
\max_{i}\{a_{i}\}$ $\leq $ $\sum_{i=1}^{p}a_{i}$. If $x_{i,t}$ is
sub-exponential then for $q$ $=$ $1$ and $\mathcal{L}_{1}$-boundedness, and
the Cauchy-Schwartz inequality, $\mathbb{E}(\max_{i}1/n%
\sum_{t=1}^{n}x_{i,t}^{2})^{1/2}$ $\leq $ $\{\mathbb{E}\max_{i,t}|x_{i,t}|%
\times n^{-1}\sum_{t=1}^{n}\mathbb{E}\max_{i}|x_{i,t}|\}^{1/2}$ $=$ $O(\sqrt{%
\ln (pn)\ln (p}))$. By restricting dependence and tail decay, other Orlicz
norms besides the one above corresponding to $x$ $\mapsto $ $x^{q}$ can be
used, including exponential under a sub-exponential assumption 
\citep[Lemma
2.2.2]{vanderVaartWellner1996}.

\citet[Lemma 8]{Chernozhukov_etal2015} use a similar symmetrization approach
to prove an $\mathcal{L}_{1}$-maximal inequality 
\begin{equation*}
\mathbb{E}\max_{i}\left\vert \bar{x}_{i,n}\right\vert \lesssim \sqrt{\frac{%
\ln (p)}{n}\times \max_{i}\frac{1}{n}\sum_{t=1}^{n}\mathbb{E}x_{i,t}^{2}}+%
\mathbb{E}\max_{i,t}x_{i,t}^{2}\frac{\ln (p)}{n}.
\end{equation*}%
See also \citet[Lemma 14.18,
Corollary 14.4]{BuhlmannVanDeGeer2011} for a low dimensional result, and see 
\cite{JuditskyNemirovski2008}, \cite{Dumbgen_etal2010} and \cite%
{MassartRossignol2013}\ for extensions for independent random variables in
Banach space. Extant results rely on square integrability of the $l_{\infty
} $ envelope $\max_{i,t}|x_{i,t}|$, and generally exploit symmetrization.

Often the goal is an optimal constant $C_{n}(p,q)$ in different spacial
contexts \citep[e.g.][]{Nemirovski2000,Dumbgen_etal2010}, typically in low
dimensions \citep[e.g.][]{IbragimovSharakhmetov1998}. Extensions to
dependent random variables have a long history in the low dimensional case;
see, e.g., \cite{vonBahrEsseen1965}, \cite{McLeish1975}, \cite{Yokoyama1980}%
, \cite{Hitczenko1990}, \cite{delaPena_etal2003}, \cite%
{Merlevede_Peligrad_2013} and \cite{Szewczak2015}. Cf. \cite{MZ1937} and 
\cite{Rosenthal1970}.

The usefulness of (\ref{Nem_q}) can be easily understood. Under $\mathcal{L}%
_{s}$-boundedness a \textquotedblleft least sharp\textquotedblright\ bound
for $\mathbb{E}\max_{i}|\bar{x}_{i,n}|^{q}$, modulo universal multiplicative
constants, is $p^{q/s}$. Strengthening the bound to $(\ln (p))^{q/2}$ under
general tail and memory conditions has great practical application for
concentration inequalities in a variety of disciplines, and can promote an
exponential bound on $p$\ yielding a max-law of large numbers $\max_{i}|\bar{%
x}_{i,n}|$ $\overset{p}{\rightarrow }$ $0$. Such LLN's appear in many places
in the high dimensional literature, including regularized estimation theory
like debiased Lasso, wavelet-like white noise tests and high dimensional
regression model tests to name just a few 
\citep[see,
e.g.,][]{Dezeure_etal_2017,Hill2024_maxlln,JinWangWang2015,HillLi2025}.

Under \textit{dependence}, however, a (classic) symmetrization argument is
apparently unavailable. We therefore approach the problem akin to the high
dimensional dependent bootstrap and Gaussian approximation literature %
\citep[e.g.][]{Mammen1993,Chernozhukov_etal2013,Chernozhukov_etal2017,Chernozhukov_etal2019,ZhangCheng2018,Dezeure_etal_2017,ZhangCheng2014,ChangJiangShao2023,ChangChenWu2024,Hill2024_maxlln,HillLi2025}%
. We initially assume $x_{i,t}$ is bounded in order to yield a version of (%
\ref{Nem_q}) by using telescoping sub-sample blocks. We then derive the main
result by approximating $\bar{x}_{i,n}$ both by negligible truncation and
blocking. We show a comparison of Gaussian distributions with and without
truncation or blocking delivers (\ref{Nem_q}), modified to account for
blocking.

Our general proof requires a new concentration bound for tail probability
measures, akin to \cite{Nemirovski2000}'s bound in $\mathcal{L}_{q}$, levied
in general and under sub-exponential and $\mathcal{L}_{q}$-boundedness
conditions. Set $\mathbb{\bar{P}}_{\mathcal{M}}$ $:=$ $\max_{i}\mathbb{P}(|%
\bar{x}_{i,n}|$ $>$ $\mathcal{M})$. In general for any $\zeta $ $\in $ $%
(0,1] $ and irrespective of dependence, we show 
\begin{equation*}
\mathbb{P}\left( \max_{i}\left\vert \bar{x}_{i,n}\right\vert >\mathcal{M}%
_{n}\right) \leq \frac{\ln (p)}{\ln (\mathbb{\bar{P}}_{\mathcal{M}%
_{n}}^{-\zeta })}+\frac{\mathbb{\bar{P}}_{\mathcal{M}_{n}}^{1-\zeta }}{\ln (%
\mathbb{\bar{P}}_{\mathcal{M}_{n}}^{-\zeta })}.
\end{equation*}%
Under $\mathcal{L}_{q}$-boundedness this becomes $\mathbb{P}(\max_{i}|\bar{x}%
_{i,n}|$ $>$ $\mathcal{M}_{n})$ $\lesssim $ $2\ln (p)/\ln (\mathcal{M}%
_{n}^{q}/\max_{i}\mathbb{E}\left\vert \bar{x}_{i,n}\right\vert ^{q})$. In a
general dependence setting this yields a significant improvement over using
Markov's inequality $\mathbb{P}(\max_{i}|\bar{x}_{i,n}|$ $>$ $\mathcal{M}%
_{n})$ $\leq $ $\mathcal{M}_{n}^{-q}\mathbb{E}\max_{i}|\bar{x}_{i,n}|^{q}$
and our main bound (\ref{b}) on $\mathbb{E}\max_{i}|\bar{x}_{i,n}|^{q}$, cf.
Proposition \ref{prop:main_unbound}. The latter entails multiple
approximation errors that are bypassed when we directly treat $\mathbb{P}%
(\max_{i}\left\vert \bar{x}_{i,n}\right\vert $ $>$ $\mathcal{M}_{n})$ as a
moment $\mathbb{E}\mathcal{I}_{\max_{i}\left\vert \bar{x}_{i,n}\right\vert >%
\mathcal{M}_{n}}$.

The main results are presented in Section \ref{sec:Lq_ineq}. Examples
involving heterogeneous geometric mixing sequences with sub-exponential
tails, and heterogeneous physical dependent sequences under $\mathcal{L}_{q}$%
-boundedness, are presented in Section \ref{sec:Exs}. Concluding remarks are
left for Section \ref{sec:conc}, and omitted proofs are relegated to the
appendix.

We assume all random variables exist on a complete measure space $(\Omega ,%
\mathcal{F},\mathbb{P})$ in order to side-step any measurability issues
concerning suprema \citep[see][Appendix C]{Pollard1984}. $\mathbb{E}$ is the
expectations operator; $\mathbb{E}_{\mathcal{A}}$ is expectations
conditional on $\mathcal{F}$-measurable $\mathcal{A}$. $\mathcal{L}_{q}$ $:=$
$\{X,$ $\sigma (X)$ $\subset $ $\mathcal{F}:$ $\mathbb{E}|X|^{q}$ $<$ $%
\infty \}$. $||\cdot ||_{q}$ is the $\mathcal{L}_{q}$-norm. $a.s.$ is $%
\mathbb{P}$-\textit{almost surely}. $K$ $>$ $0$ is a finite constant that
may have different values in different places. $x$ $\lesssim $ $y$ if $x$ $%
\leq $ $Ky$ for some $K$ $>$ $0$ that is not a function of $n,x,y$. $%
\mathcal{I}_{A}$ is the indicator function: $\mathcal{I}_{A}$ $=$ $1$ if $A$
is true. $|z|_{+}$ $=$ $0\vee z$. $[\cdot ]_{+}$ rounds to the next greater
integer.

\section{$\mathcal{L}_{q}$-maximal inequality\label{sec:Lq_ineq}}

It is useful to recall how symmetrization works under independence. Write $%
\mathfrak{X}^{(n)}$ $:=$ $\{x_{t}\}_{t=1}^{n}$, and let $\{x_{t}^{\ast }\}$
be an independent copy of $\{x_{t}\}$. Then $\mathbb{E}x_{i,t}^{\ast }$ $=$ $%
\mathbb{E}x_{i,t}$ $=$ $0$, hence by the conditional Jensen's inequality and
Fubini's theorem%
\begin{equation*}
\left\Vert \max_{i}\left\vert \bar{x}_{i,n}\right\vert \right\Vert _{q}\leq
\left\Vert \mathbb{E}_{\mathfrak{X}^{(n)}}\max_{i}\left\vert \frac{1}{n}%
\sum_{t=1}^{n}\left( x_{i,t}-x_{i,t}^{\ast }\right) \right\vert \right\Vert
_{q}\leq \left\Vert \max_{i}\left\vert \frac{1}{n}\sum_{t=1}^{n}\left(
x_{i,t}-x_{i,t}^{\ast }\right) \right\vert \right\Vert _{q}.
\end{equation*}%
If $\{\varepsilon _{t}\}$ are iid Rademacher (i.e. $\mathbb{P}(\varepsilon
_{t}$ $=$ $-1)$ $=$ $\mathbb{P}(\varepsilon _{t}$ $=$ $1)$ $=$ $1/2$),
independent of $\{x_{t}\}$, then symmetrically distributed $%
x_{i,t}-x_{i,t}^{\ast }$ and $\varepsilon _{t}\left( x_{i,t}-x_{i,t}^{\ast
}\right) $ have the same distribution. Indeed, under mutual and serial
independence $\max_{i}|1/n\sum_{t=1}^{n}(x_{i,t}$ $-$ $x_{i,t}^{\ast })|$
and $\max_{i}|1/n\sum_{t=1}^{n}\varepsilon _{t}(x_{i,t}$ $-$ $x_{i,t}^{\ast
})|$ have the same distribution. Thus $||\max_{i}|1/n\sum_{t=1}^{n}(x_{i,t}$ 
$-$ $x_{i,t}^{\ast })|||_{q}$ $=$ $||\max_{i}|1/n\sum_{t=1}^{n}\varepsilon
_{t}(x_{i,t}$ $-$ $x_{i,t}^{\ast })|||_{q}$. The latter is easily bounded
using triangle and Minkowski inequalities yielding the \textquotedblleft
desymmetrized\textquotedblright\ $||\max_{i}|1/n\sum_{t=1}^{n}\varepsilon
_{t}(x_{i,t}$ $-$ $x_{i,t}^{\ast })|||_{q}$ $\leq $ $2||\max_{i}|\bar{x}%
_{i,n}|||_{q}$. Now note that $\varepsilon _{t}(x_{i,t}$ $-$ $x_{i,t}^{\ast
})$ is sub-Gaussian once we condition on $\{x_{i,t},x_{i,t}^{\ast
}\}_{t=1}^{n}$. A classic log-exp bound and Hoeffding's inequality are then
enough to yield (\ref{Nem_q}), cf. 
\citet[Lemmas 14.10, 14.14,
14.24]{BuhlmannVanDeGeer2011}.

Under \textit{dependence} we use telescoping blocks and cross-block
independent multipliers. Let $b_{n}$ $\in $ $\{1,...,n$ $-$ $1\}$ be a
pre-set block size, $b_{n}$ $\rightarrow $ $\infty $, $b_{n}$ $=$ $o(n)$.
Define the number of blocks $\mathcal{N}_{n}$ $:=$ $[n/b_{n}]$, and index
sets $\mathfrak{B}_{l}$ $:=$ $\{(l$ $-$ $1)b_{n}$ $+$ $1,\dots ,lb_{n}\}$
with $l$ $=$ $1,\dots ,\mathcal{N}_{n}$. Assume $\mathcal{N}_{n}b_{n}$ $=$ $%
n $ throughout to reduce notation. Generate bounded iid random variables $%
\{\varepsilon _{l}\}_{l=1}^{\mathcal{N}_{n}}$ independent of $%
\{x_{t}\}_{t=1}^{n}$, where $\mathbb{P}(|\varepsilon _{l}|$ $\leq $ $c)$ $=$ 
$1$, $c$ $\in $ $(0,\infty )$, and define a sequence of multipliers $\{\eta
_{t}\}_{t=1}^{n}$ by $\eta _{t}$ $=$ $\varepsilon _{l}$ if $t$ $\in $ $%
\mathfrak{B}_{l}$. Define partial sums and a variance%
\begin{eqnarray}
&&\mathcal{X}_{n}(i):=\frac{1}{\sqrt{n}}\sum_{t=1}^{n}x_{i,t}\text{, \ }%
\sigma _{n}^{2}(i,j):=\mathbb{E}\mathcal{X}_{n}(i)\mathcal{X}_{n}(j)\text{
and }\sigma _{n}^{2}(i):=\mathbb{E}\mathcal{X}_{n}^{2}(i)  \label{X_sig} \\
&&\mathcal{X}_{n}^{\ast }(i):=\frac{1}{\sqrt{n}}\sum_{t=1}^{n}\eta
_{t}x_{i,t}=\frac{1}{\sqrt{n}}\sum_{l=1}^{\mathcal{N}_{n}}\varepsilon _{l}%
\mathcal{S}_{n,l}(i)\text{ where }\mathcal{S}_{n,l}(i):=%
\sum_{t=(l-1)b_{n}+1}^{lb_{n}}x_{i,t}.\text{ \ \ \ \ \ \ \ }  \notag
\end{eqnarray}%
Thus $\mathcal{X}_{n}^{\ast }(i)$ is a multiplier (wild) bootstrap
approximation of $\mathcal{X}_{n}(i)$. It is, very loosely speaking, a
\textquotedblleft symmetrized\textquotedblright\ version of $\mathcal{X}%
_{n}(i)$. As long as $b_{n}$ $\rightarrow $ $\infty $ and $b_{n}$ $=$ $o(n)$
then under a broad array of dependence settings $\mathcal{X}_{n}^{\ast }(i)$
is a draw from the distribution governing $\mathcal{X}_{n}(i)$,
asymptotically with probability approaching one %
\citep[e.g.][]{Liu1988,GineZinn1990,HillLi2025}. A Rademacher assumption for 
$\varepsilon _{l}$ under dependence cannot precisely serve the purpose it
does under independence. We still require boundedness, however, in order to
invoke Hoeffding's inequality in applications.

Now let $\{\boldsymbol{X}_{n}(i)\}_{i=1}^{p}$ be a Gaussian process, $%
\boldsymbol{X}_{n}(i)$ $\sim $ $N(0,\sigma _{n}^{2}(i))$, and define%
\begin{eqnarray}
&&\rho _{n}=\rho (n,p):=\sup_{z\geq 0}\left\vert \mathbb{P}\left(
\max_{i}\left\vert \mathcal{X}_{n}(i)\right\vert \leq z\right) -\mathbb{P}%
\left( \max_{i}\left\vert \boldsymbol{X}_{n}(i)\right\vert \leq z\right)
\right\vert  \label{rho_n_dfn} \\
&&\rho _{n}^{\ast }=\rho ^{\ast }(n,p):=\sup_{z\geq 0}\left\vert \mathbb{P}%
\left( \max_{i}\left\vert \mathcal{X}_{n}^{\ast }(i)\right\vert \leq
z\right) -\mathbb{P}\left( \max_{i}\left\vert \boldsymbol{X}%
_{n}(i)\right\vert \leq z\right) \right\vert .  \notag
\end{eqnarray}%
Hence $\sup_{z\geq 0}|\mathbb{P}(\max_{i}|\mathcal{X}_{n}^{\ast }(i)|$ $\leq 
$ $z)$ $-$ $\mathbb{P}(\max_{i}|\mathcal{X}_{n}(i)|$ $\leq $ $z)|$ $\leq $ $%
\rho _{n}$ $+$ $\rho _{n}^{\ast }$ by the triangle inequality: blocking has
an asymptotically negligible (in probability) impact whenever the Kolmogorov
distances $\{\rho _{n},\rho _{n}^{\ast }\}$ $\rightarrow $ $0$. The latter
holds in many settings: see Examples \ref{sec:Ex1}-\ref{sec:Ex3}, and see,
e.g., \cite{Chernozhukov_etal2013,Chernozhukov_etal2017}, \cite%
{JinWangWang2015}, \cite{ZhangWu2017}, \cite{ZhangCheng2018}, \cite%
{ChangJiangShao2023}, and \cite{ChangChenWu2024}.

\subsection{Bounded\label{sec:bound}}

We have our first \cite{Nemirovski2000}-like moment bound under boundedness
and arbitrary dependence. The result provides a foundation for using a
negligible truncation approximation under unboundedness.

\begin{proposition}
\label{prop:main} Let $\mathbb{P}(\max_{i,t}\left\vert x_{i,t}\right\vert $ $%
\leq $ $\mathcal{M}_{n})$ $=$ $1$ for some non-decreasing sequence $\{%
\mathcal{M}_{n}\}$, $\mathcal{M}_{n}$ $\in $ $(0,\infty )$ where $\mathcal{M}%
_{n}$ $\rightarrow $ $\infty $ is possible. Then%
\begin{eqnarray}
\mathbb{E}\max_{i}\left\vert \bar{x}_{i,n}\right\vert ^{q} &\leq &\left( 
\frac{2c^{2}\ln \left( 2p\right) }{\mathcal{N}_{n}}\right) ^{q/2}\sqrt{%
\mathbb{E}\max_{i,l}\left\vert \frac{1}{b_{n}}\mathcal{S}_{n,l}(i)\right%
\vert ^{q}\mathbb{E}\left( \max_{i}\frac{1}{\mathcal{N}_{n}}\sum_{l=1}^{%
\mathcal{N}_{n}}\left\vert \mathcal{S}_{n,l}(i)\right\vert \right) ^{q}}
\label{main} \\
&&\text{ \ \ \ }+\left( \frac{\mathcal{M}_{n}}{\sqrt{n}}\right) ^{q}\left\{
\rho _{n}+\rho _{n}^{\ast }\right\} .  \notag
\end{eqnarray}
\end{proposition}

\begin{remark}
\normalfont Use Minkowski and Jensen inequalities to yield simple refinements%
\begin{eqnarray*}
\mathbb{E}\max_{i}\left\vert \bar{x}_{i,n}\right\vert ^{q} &\leq &\left( 
\frac{2c^{2}\ln \left( 2p\right) }{\mathcal{N}_{n}}\right) ^{q/2}\sqrt{%
\mathbb{E}\max_{i,t}\left\vert x_{i,t}\right\vert ^{q}\times \mathbb{E}%
\max_{i,l}\left( \frac{1}{b_{n}}\sum_{t=(l-1)b_{n}+1}^{lb_{n}}\left\vert
x_{i,t}\right\vert \right) ^{q}}\text{ \ \ \ } \\
&&\text{ \ \ \ }+\left( \frac{\mathcal{M}_{n}}{\sqrt{n}}\right) ^{q}\left\{
\rho _{n}+\rho _{n}^{\ast }\right\} \\
&& \\
&\leq &\left( \frac{2c^{2}\ln \left( 2p\right) }{\mathcal{N}_{n}}\right)
^{q/2}\mathbb{E}\max_{i,t}\left\vert x_{i,t}\right\vert ^{q}+\left( \frac{%
\mathcal{M}_{n}}{\sqrt{n}}\right) ^{q}\left\{ \rho _{n}+\rho _{n}^{\ast
}\right\} .
\end{eqnarray*}%
If $\varepsilon _{l}$ is Rademacher then $c$ $=$ $1$, while in general
applications we need $c$ $\geq $ $1$. Indeed, we are not free to choose $c$
and thus $\varepsilon _{l}$. In order to verify $\{\rho _{n},\rho _{n}^{\ast
}\}$ $\rightarrow $ $0$ we employ a Gaussian-to-Gaussian comparison that
exploits $\mathbb{E}\varepsilon _{l}^{2}$ $=$ $1$, a natural requirement in
the multiplier bootstrap literature, hence $c$ $\geq $ $1$.
\end{remark}

\begin{remark}
\normalfont If $q$ $=$ $1$ then we only need $\mathbb{E}\max_{i,t}|x_{i,t}|$ 
$<$ $\infty $, although with boundedness this is immaterial since all
moments exist. It becomes important, however, when we relax boundedness
where a second moment need not exist.
\end{remark}

\begin{remark}
\label{rm:Kol_M}\normalfont Boundedness is used solely to yield a bounded
integral under general dependence $\mathbb{E}\max_{i}\left\vert \bar{x}%
_{i,n}\right\vert ^{q}$ $=$ $\int_{0}^{\mathcal{M}_{n}}\mathbb{P}%
(\max_{i}|\max_{i}|\bar{x}_{i,n}||^{q}$ $>$ $u)du$. The Kolmogorov distances 
$\{\rho _{n},\rho _{n}^{\ast }\}$ themselves, however, do not require nor
exploit boundedness.
\end{remark}

\noindent \textbf{Proof}. By the triangle inequality and definitions of $%
\{\rho _{n},\rho _{n}^{\ast }\},$%
\begin{equation*}
\sup_{z\geq 0}\left\vert \mathbb{P}\left( \max_{i}\left\vert \frac{1}{\sqrt{n%
}}\sum_{t=1}^{n}x_{i,t}\right\vert \leq z\right) -\mathbb{P}\left(
\max_{i}\left\vert \frac{1}{\sqrt{n}}\sum_{l=1}^{\mathcal{N}_{n}}\varepsilon
_{l}\mathcal{S}_{n,l}(i)\right\vert \leq z\right) \right\vert \leq \rho
_{n}+\rho _{n}^{\ast }.
\end{equation*}%
Now replace $x_{i,t}$ with $x_{i,t}/\sqrt{n}$ to yield for each $x$ $\geq $ $%
0$, 
\begin{eqnarray}
&&\left\vert \mathbb{P}\left( \max_{i}\left\vert \bar{x}_{i,n}\right\vert
\leq x\right) -\mathbb{P}\left( \max_{i}\left\vert \frac{1}{n}\sum_{l=1}^{%
\mathcal{N}_{n}}\varepsilon _{l}\mathcal{S}_{n,l}(i)\right\vert \leq
x\right) \right\vert  \label{ppgg} \\
&&\text{ \ \ \ \ \ }=\left\vert \mathbb{P}\left( \max_{i}\left\vert
\sum_{t=1}^{n}\frac{x_{i,t}}{\sqrt{n}}\right\vert \leq \sqrt{n}x\right) -%
\mathbb{P}\left( \max_{i}\left\vert \sum_{l=1}^{\mathcal{N}_{n}}\varepsilon
_{l}\frac{\mathcal{S}_{n,l}(i)}{\sqrt{n}}\right\vert \leq \sqrt{n}x\right)
\right\vert \leq \rho _{n}+\rho _{n}^{\ast }\text{.}  \notag
\end{eqnarray}%
Recall $\mathfrak{X}^{(n)}$ $:=$ $\{x_{t}\}_{t=1}^{n}$. Boundedness, twice a
change of variables, and (\ref{ppgg}) imply%
\begin{eqnarray}
\mathbb{E}\max_{i}\left\vert \bar{x}_{i,n}\right\vert ^{q} &=&\int_{0}^{%
\mathcal{M}_{n}}\mathbb{P}\left( \max_{i}\left\vert \max_{i}\left\vert \bar{x%
}_{i,n}\right\vert \right\vert ^{q}>u\right) du  \notag \\
&=&q\int_{0}^{\mathcal{M}_{n}}v^{q-1}\mathbb{P}\left( \max_{i}\left\vert 
\sqrt{n}\bar{x}_{i,n}\right\vert >\sqrt{n}v\right) dv  \notag \\
&=&\frac{q}{n^{q/2}}\int_{0}^{\sqrt{n}\mathcal{M}_{n}}z^{q-1}\mathbb{P}%
\left( \max_{i}\left\vert \sqrt{n}\bar{x}_{i,n}\right\vert >z\right) dz 
\notag \\
&\leq &\frac{q}{n^{q/2}}\int_{0}^{\sqrt{n}\mathcal{M}_{n}}z^{q-1}\mathbb{P}%
\left( \max_{i}\left\vert \frac{1}{\sqrt{n}}\sum_{l=1}^{\mathcal{N}%
_{n}}\varepsilon _{l}\mathcal{S}_{n,l}(i)\right\vert \leq z\right) dz  \notag
\\
&&\text{ \ \ \ }+\frac{\rho _{n}\vee \rho _{n}^{\ast }}{n^{q/2}}q\int_{0}^{%
\mathcal{M}_{n}}u^{q-1}du  \notag \\
&&  \notag \\
&=&\mathbb{E}\left( \mathbb{E}_{\mathfrak{X}^{(n)}}\max_{i}\left\vert \frac{1%
}{n}\sum_{l=1}^{\mathcal{N}_{n}}\varepsilon _{l}\mathcal{S}%
_{n,l}(i)\right\vert ^{q}\right) +\left( \frac{\mathcal{M}_{n}}{\sqrt{n}}%
\right) ^{q}\left\{ \rho _{n}+\rho _{n}^{\ast }\right\} .  \label{EEx}
\end{eqnarray}%
Since $\varepsilon _{l}$ is iid and $\mathbb{P}(|\varepsilon _{l}|$ $\leq $ $%
c)$ $=$ $1$, $c$ $\in $ $(0,\infty )$, Hoeffding's inequality applies to $%
\mathbb{E}_{\mathfrak{X}^{(n)}}(\cdot )$ with a log-exp bound to yield 
\citep[see][Lemma
14.14]{BuhlmannVanDeGeer2011} 
\begin{equation*}
\mathbb{E}_{\mathfrak{X}^{(n)}}\max_{i}\left\vert \frac{1}{n}\sum_{l=1}^{%
\mathcal{N}_{n}}\varepsilon _{l}\mathcal{S}_{n,l}(i)\right\vert ^{q}\leq
\left( \frac{2c^{2}\ln \left( 2p\right) }{n}\right) ^{q/2}\left( \max_{i}%
\frac{1}{n}\sum_{l=1}^{\mathcal{N}_{n}}\mathcal{S}_{n,l}^{2}(i)\right)
^{q/2}.
\end{equation*}%
See, e.g., \cite{Bentkus2004,Bentkus2008} for Hoeffding-related bounds
allowing for unbounded $\varepsilon _{l}$. Hence by Fubini's theorem, the
Cauchy-Schwartz inequality and $n$ $=$ $\mathcal{N}_{n}b_{n}$,%
\begin{eqnarray}
&&\text{ \ \ \ }\mathbb{E}\max_{i}\left\vert \frac{1}{n}\sum_{l=1}^{\mathcal{%
N}_{n}}\varepsilon _{l}\mathcal{S}_{n,l}(i)\right\vert ^{q}  \label{Esum*1}
\\
&&\text{ \ \ \ \ \ \ \ \ \ \ \ }\leq \left( \frac{2c^{2}\ln \left( 2p\right) 
}{\mathcal{N}_{n}}\right) ^{q/2}\mathbb{E}\left( \max_{i}\frac{1}{n}%
\sum_{l=1}^{\mathcal{N}_{n}}\mathcal{S}_{n,l}^{2}(i)\right) ^{q/2}  \notag \\
&&\text{ \ \ \ \ \ \ \ \ \ \ \ }\leq \left( \frac{2c^{2}\ln \left( 2p\right) 
}{\mathcal{N}_{n}}\right) ^{q/2}\mathbb{E}\left[ \max_{i,l}\left\vert \frac{1%
}{b_{n}}\mathcal{S}_{n,l}(i)\right\vert ^{q/2}\left( \max_{i}\frac{1}{%
\mathcal{N}_{n}}\sum_{l=1}^{\mathcal{N}_{n}}\left\vert \mathcal{S}%
_{n,l}(i)\right\vert \right) ^{q/2}\right]  \notag \\
&&\text{ \ \ \ \ \ \ \ \ \ \ \ }\leq \left( \frac{2c^{2}\ln \left( 2p\right) 
}{\mathcal{N}_{n}}\right) ^{q/2}\sqrt{\mathbb{E}\max_{i,l}\left\vert \frac{1%
}{b_{n}}\mathcal{S}_{n,l}(i)\right\vert ^{q}\mathbb{E}\left( \max_{i}\frac{1%
}{\mathcal{N}_{n}}\sum_{l=1}^{\mathcal{N}_{n}}\left\vert \mathcal{S}%
_{n,l}(i)\right\vert \right) ^{q}}.  \notag
\end{eqnarray}%
Combine (\ref{EEx}) with (\ref{Esum*1}) to conclude (\ref{main}). $\mathcal{%
QED}$.

\subsection{Unbounded\label{sec:unbound}}

Now assume $x_{t}$\ has support $\mathbb{R}^{p}$, and let $\mathcal{\{M}%
_{n}\}_{n\in \mathbb{N}}$ be a sequence of monotonically increasing positive
reals, $\mathcal{M}_{n}$ $\rightarrow $ $\infty $. We approximate $x_{t}$
with a centered negligibly truncated version 
\begin{equation*}
y_{i,t}^{(\mathcal{M})}:=x_{i,t}^{(\mathcal{M})}-\mathbb{E}x_{i,t}^{(%
\mathcal{M})}\text{ \ where }x_{i,t}^{(\mathcal{M})}:=x_{i,t}\mathcal{I}%
_{\left\vert x_{i,t}\right\vert \leq \mathcal{M}_{n}}+\mathcal{M}_{n}%
\mathcal{I}_{\left\vert x_{i,t}\right\vert >\mathcal{M}_{n}}\text{.}
\end{equation*}%
Define partial sums under truncation with and without blocking, and a
partial sum (co)variance%
\begin{eqnarray*}
&&\mathcal{X}_{n}^{(\mathcal{M})}(i):=\frac{1}{\sqrt{n}}%
\sum_{t=1}^{n}y_{i,t}^{(\mathcal{M})}\text{, \ }\sigma _{n}^{(\mathcal{M}%
)2}(i,j):=\mathbb{E}\mathcal{X}_{n}^{(\mathcal{M})}(i)\mathcal{X}_{n}^{(%
\mathcal{M})}(j)\text{ and }\sigma _{n}^{(\mathcal{M})2}(i)=\sigma _{n}^{(%
\mathcal{M})2}(i,i) \\
&&\mathcal{X}_{n}^{(\mathcal{M})\ast }(i):=\frac{1}{\sqrt{n}}%
\sum_{t=1}^{n}\eta _{t}y_{i,t}^{(\mathcal{M})}=\frac{1}{\sqrt{n}}\sum_{l=1}^{%
\mathcal{N}_{n}}\varepsilon _{l}\mathcal{S}_{n,l}^{(\mathcal{M})}(i)\text{
where }\mathcal{S}_{n,l}^{(\mathcal{M})}(i):=%
\sum_{t=(l-1)b_{n}+1}^{lb_{n}}y_{i,t}^{(\mathcal{M})}.
\end{eqnarray*}%
Let $\{\boldsymbol{X}_{n}^{(\mathcal{M})}(i)\}_{i\in \mathbb{N}}$ be a
Gaussian process, $\boldsymbol{X}_{n}^{(\mathcal{M})}(i)$ $\sim $ $%
N(0,\sigma _{n}^{(\mathcal{M})2}(i))$, recall $\boldsymbol{X}_{n}(i)$ $\sim $
$N(0,\mathbb{E}\mathcal{X}_{n}^{2}(i))$, and define 
\begin{eqnarray*}
&&\rho _{n}^{\ast (\mathcal{M})}:=\sup_{z\geq 0}\left\vert \mathbb{P}\left(
\max_{i}\left\vert \mathcal{X}_{n}^{(\mathcal{M})\ast }(i)\right\vert \leq
z\right) -\mathbb{P}\left( \max_{i}\left\vert \boldsymbol{X}_{n}^{(\mathcal{M%
})}(i)\right\vert \leq z\right) \right\vert \\
&&\delta _{n}^{(\mathcal{M})}:=\sup_{z\geq 0}\left\vert \mathbb{P}\left(
\max_{i}\left\vert \boldsymbol{X}_{n}^{(\mathcal{M})}(i)\right\vert \leq
z\right) -\mathbb{P}\left( \max_{i}\left\vert \boldsymbol{X}%
_{n}(i)\right\vert \leq z\right) \right\vert .
\end{eqnarray*}%
Thus $\rho _{n}^{\ast (\mathcal{M})}$ measure how close the truncated and
blocked $\max_{i}|\mathcal{X}_{n}^{(\mathcal{M})\ast }(i)|$ is to the
max-Gaussian law under truncation $\max_{i}|\boldsymbol{X}_{n}^{(\mathcal{M}%
)}(i)|$, and $\delta _{n}^{(\mathcal{M})}$ captures the distance between
Gaussian laws with and without truncation. Notice inequality (\ref{main})
instantly holds for the bounded $\mathcal{X}_{n}^{(\mathcal{M})}(i)$ by
Proposition \ref{prop:main}.

Recall $\sigma _{n}^{2}(i,j)$ $:=$ $\mathbb{E}\mathcal{X}_{n}(i)\mathcal{X}%
_{n}(j)$. By Lemma C.5 in \cite{Chen_Kato_2019}, cf. Theorem 2 in \cite%
{Chernozhukov_etal2015},

\begin{equation*}
\delta _{n}^{(\mathcal{M})}\lesssim \Delta _{n}^{(\mathcal{M})1/3}\times \ln
(p)^{2/3}\text{ \ where \ }\Delta _{n}^{(\mathcal{M})}:=\max_{i,j}\left\vert
\sigma _{n}^{2}(i,j)-\sigma _{n}^{(\mathcal{M})2}(i,j)\right\vert .
\end{equation*}%
It is straightforward to show $\delta _{n}^{(\mathcal{M})}$ $\rightarrow $ $%
0 $ by manipulating the truncation points $\{\mathcal{M}_{n}\}$ and bounding 
$\ln (p)$, without any reference to underlying dependence or tail decay
other than the existence of a higher moment. First, by construction and the
triangle inequality%
\begin{eqnarray*}
&&\Delta _{n}^{(\mathcal{M})}\leq \max_{i,j}\left\vert \frac{1}{n}%
\sum_{s,t=1}^{n}\mathbb{E}x_{i,s}x_{i,t}-\frac{1}{n}\sum_{s,t=1}^{n}\mathbb{E%
}x_{i,s}\mathcal{I}_{\left\vert x_{i,s}\right\vert \leq \mathcal{M}%
_{n}}x_{i,t}\mathcal{I}_{\left\vert x_{i,t}\right\vert \leq \mathcal{M}%
_{n}}\right\vert \\
&&\text{ \ \ \ \ \ \ \ \ \ \ }+2\mathcal{M}_{n}\max_{i,j}\left\vert \frac{1}{%
n}\sum_{s,t=1}^{n}\mathbb{E}\mathcal{I}_{\left\vert x_{i,s}\right\vert >%
\mathcal{M}_{n}}\left( x_{i,t}\mathcal{I}_{\left\vert x_{i,t}\right\vert
\leq \mathcal{M}_{n}}+\mathcal{M}_{n}\mathcal{I}_{\left\vert
x_{i,t}\right\vert >\mathcal{M}_{n}}\right) \right\vert :=\mathfrak{A}_{n}+%
\mathfrak{B}_{n}.
\end{eqnarray*}%
Now suppose $\mathbb{E}|x_{i,t}|^{2r}$ $<$ $\infty $ for some $r$ $>$ $1$,
and assume $\mathcal{M}_{n}$ $\rightarrow $ $\infty $ sufficiently fast, 
\begin{equation}
\frac{n^{1/(r-1)}\left\{ 1+\max_{i,t}\left\Vert x_{i,s}\right\Vert
_{2r}\right\} ^{2r/(r-1)}}{\mathcal{M}_{n}}\rightarrow 0.  \label{Mn}
\end{equation}%
Thus under data homogeneity $\mathcal{M}_{n}/n^{1/(r-1)}$ $\rightarrow $ $%
\infty $. By multiple applications of Cauchy-Schwartz, H\"{o}lder, Lyapunov
and Markov inequalities, and (\ref{Mn}) with $r$ $>$ $1$, 
\begin{eqnarray*}
\mathfrak{A}_{n} &\leq &\left\vert \frac{1}{n}\sum_{s,t=1}^{n}\mathbb{E}%
x_{i,t}x_{i,s}\mathcal{I}_{\left\vert x_{i,s}\right\vert >\mathcal{M}%
_{n}}\right\vert +\left\vert \frac{1}{n}\sum_{s,t=1}^{n}\mathbb{E}x_{i,s}%
\mathcal{I}_{\left\vert x_{i,s}\right\vert \leq \mathcal{M}_{n}}x_{i,t}%
\mathcal{I}_{\left\vert x_{i,t}\right\vert >\mathcal{M}_{n}}\right\vert \\
&\leq &\frac{1}{n}\sum_{s,t=1}^{n}\left( \mathbb{E}x_{i,t}^{2}\right)
^{1/2}\left( \mathbb{E}x_{i,s}^{2}\mathcal{I}_{\left\vert x_{i,s}\right\vert
>\mathcal{M}_{n}}\right) ^{1/2}+\frac{1}{n}\sum_{s,t=1}^{n}\left( \mathbb{E}%
x_{i,s}^{2}\mathcal{I}_{\left\vert x_{i,s}\right\vert \leq \mathcal{M}%
_{n}}\right) ^{1/2}\left( \mathbb{E}x_{i,t}^{2}\mathcal{I}_{\left\vert
x_{i,t}\right\vert >\mathcal{M}_{n}}\right) ^{1/2} \\
&\leq &\frac{1}{n}\sum_{s,t=1}^{n}\left( \mathbb{E}x_{i,t}^{2}\right)
^{1/2}\left( \mathbb{E}x_{i,s}^{2r}\right) ^{1/(2r)}\mathbb{P}\left(
\left\vert x_{i,s}\right\vert >\mathcal{M}_{n}\right) ^{(r-1)/(r2)} \\
&&\text{ \ \ \ }+\frac{1}{n}\sum_{s,t=1}^{n}\left( \mathbb{E}x_{i,s}^{2}%
\mathcal{I}_{\left\vert x_{i,s}\right\vert \leq \mathcal{M}_{n}}\right)
^{1/2}\left( \mathbb{E}x_{i,s}^{2r}\right) ^{1/(2r)}\mathbb{P}\left(
\left\vert x_{i,s}\right\vert >\mathcal{M}_{n}\right) ^{(r-1)/(r2)} \\
&& \\
&\leq &n\mathcal{M}_{n}^{-(r-1)}\max_{i,t}\left( \mathbb{E}%
x_{i,t}^{2}\right) ^{1/2}\max_{i,t}\left( \mathbb{E}x_{i,t}^{2r}\right)
^{1/(2r)}\max_{i,t}\left( \mathbb{E}\left\vert x_{i,t}\right\vert
^{2r}\right) ^{(r-1)/(r2)} \\
&\leq &\frac{n\max_{i,t}\left\Vert x_{i,t}\right\Vert _{2r}^{r+1}}{\mathcal{M%
}_{n}^{r-1}}\rightarrow 0.
\end{eqnarray*}%
Similarly, for some $s$ $:=$ $2r$ $>$ $2$ use (\ref{Mn}) to yield%
\begin{eqnarray*}
\mathfrak{B}_{n} &\simeq &\mathcal{M}_{n}\left\vert \frac{1}{n}%
\sum_{s,t=1}^{n}\mathbb{E}\mathcal{M}_{n}\mathcal{I}_{\left\vert
x_{i,s}\right\vert >\mathcal{M}_{n}}\left( x_{i,t}\mathcal{I}_{\left\vert
x_{i,t}\right\vert \leq \mathcal{M}_{n}}+\mathcal{M}_{n}\mathcal{I}%
_{\left\vert x_{i,t}\right\vert >\mathcal{M}_{n}}\right) \right\vert \\
&\leq &\mathcal{M}_{n}\left\vert \frac{1}{n}\sum_{s,t=1}^{n}\left( \mathbb{E}%
\left\vert x_{i,t}\right\vert ^{s}\right) ^{1/s}\mathbb{P}\left( \left\vert
x_{i,s}\right\vert >\mathcal{M}_{n}\right) ^{(s-1)/s}\right\vert \\
&\leq &\mathcal{M}_{n}\left\vert \frac{1}{n}\sum_{s,t=1}^{n}\left( \mathbb{E}%
\left\vert x_{i,t}\right\vert ^{s}\right) ^{1/s}\left( \mathbb{E}\left\vert
x_{i,s}\right\vert ^{s}\right) ^{(s-1)/s}\frac{1}{\mathcal{M}_{n}^{s-1}}%
\right\vert \\
&\leq &\frac{n}{\mathcal{M}_{n}^{s-2}}\max_{i,t}\left( \mathbb{E}\left\vert
x_{i,s}\right\vert ^{s}\right) =\frac{n}{\mathcal{M}_{n}^{2(r-1)}}\max_{i,t}%
\mathbb{E}\left\vert x_{i,s}\right\vert ^{2r}\rightarrow 0.
\end{eqnarray*}%
Thus $\Delta _{n}^{(\mathcal{M})}$ $\leq $ $\mathfrak{A}_{n}$ $+$ $\mathfrak{%
B}_{n}$ $\rightarrow $ $0$ under (\ref{Mn}). Hence when $\ln (p)$ $=$ $o(%
\mathcal{M}_{n}^{(r-1)/2}/[\sqrt{n}\max_{i,t}||x_{i,t}||_{2r}^{r}])$, for
some $r$ $>$ $1$ 
\begin{equation}
\delta _{n}^{(\mathcal{M})}\lesssim \left\{ \frac{n^{1/2}}{\mathcal{M}%
_{n}^{(r-1)/2}}\max_{i,t}\left\Vert x_{i,t}\right\Vert _{2r}^{r}\ln
(p)\right\} ^{2/3}\rightarrow 0.  \label{dM_U}
\end{equation}%
If marginal distribution tail information for $x_{i,t}$ is available then $%
\mathbb{P}(|x_{i,t}|$ $>$ $\mathcal{M}_{n})$ can be sharpened. We could
likewise exploit dependence decay, if available, to sharpen bounds on $%
1/n\sum_{s,t=1}^{n}\mathbb{E}x_{i,t}x_{i,s}\mathcal{I}_{\left\vert
x_{i,s}\right\vert >\mathcal{M}_{n}}$ and $1/n\sum_{s,t=1}^{n}\mathbb{E}%
x_{i,s}\mathcal{I}_{\left\vert x_{i,s}\right\vert \leq \mathcal{M}%
_{n}}x_{i,t}\mathcal{I}_{\left\vert x_{i,t}\right\vert >\mathcal{M}_{n}}$.

\begin{assumption}
\label{assum:unbound}Let $q$ $\geq $ $1$.$\medskip $\newline
$(a)$ $\mathbb{E}\max_{i,t}|x_{i,t}|^{(q\vee 2)r}$ $<$ $\infty $ for some $r$
$>$ $1$ and each $n$.$\medskip $\newline
$(b)$ Truncation points $\{\mathcal{M}_{n}\}$ satisfy (\ref{Mn}).
\end{assumption}

\begin{remark}
\normalfont We need $\mathbb{E}\max_{i,t}|x_{i,t}|^{(q\vee 2)r}$ $<$ $\infty 
$ for some $r$ $>$ $1$ in order for the truncation error to vanish rapidly
enough in high dimension.
\end{remark}

\begin{remark}
\normalfont(\ref{Mn}) ensures a truncation and non-truncation Gaussian
comparison \linebreak $\max_{i,j}|\sigma _{n}^{2}(i,j)$ $-$ $\sigma _{n}^{(%
\mathcal{M})2}(i,j)|$ $\rightarrow $ $0$ holds. Since symmetrization cannot
be used, it ensures remainder terms $\{\mathcal{B}_{n}(p),\mathcal{D}%
_{n}(p)\}$ defined in (\ref{BnCn}) and (\ref{Dn_gen})-(\ref{Dn_sube}) below
are negligible when $p$ is bounded. A third remainder $\mathcal{C}_{n}(p)$ $%
:=$ $\rho _{n}$ $+$ $\rho _{n}^{\ast (\mathcal{M})}$ in (\ref{BnCn}) is
negligible once a dependence setting is imposed, as in Examples \ref{sec:Ex2}
and \ref{sec:Ex3}.
\end{remark}

Define $\mathbb{\bar{P}}_{\mathcal{M}}$ $:=$ $\max_{i}\mathbb{P}(|\bar{x}%
_{i,n}|$ $>$ $\mathcal{M})$ and for $r$ $>$ $1$%
\begin{equation}
\mathcal{B}_{n}(p):=\left\{ \frac{\max_{i,t}\left\Vert x_{i,t}\right\Vert
_{2r}^{r}}{\mathcal{M}_{n}^{(r-1)}}\ln (p)\right\} ^{2/3}\text{ \ and \ }%
\mathcal{C}_{n}(p):=\rho _{n}+\rho _{n}^{\ast (\mathcal{M})}.  \label{BnCn}
\end{equation}%
We now have the main result of this paper.

\begin{proposition}
\label{prop:main_unbound} For any random variables $\{x_{t}\}_{t=1}^{n}$
satisfying Assumption \ref{assum:unbound},%
\begin{eqnarray}
&&\mathbb{E}\max_{i}\left\vert \bar{x}_{i,n}\right\vert ^{q}\lesssim \left( 
\frac{2\ln \left( 2p\right) }{\mathcal{N}_{n}}\right) ^{q/2}\sqrt{\mathbb{E}%
\max_{i,l}\left\vert \frac{1}{b_{n}}\mathcal{S}_{n,l}^{(\mathcal{M}%
)}(i)\right\vert ^{q}\mathbb{E}\left( \max_{i}\frac{1}{\mathcal{N}_{n}}%
\sum_{l=1}^{\mathcal{N}_{n}}\left\vert \mathcal{S}_{n,l}^{(\mathcal{M}%
)}(i)\right\vert \right) ^{q}}\text{ \ \ \ \ \ \ \ }  \label{b} \\
&&\text{ \ \ \ \ \ \ \ \ \ \ \ \ \ \ \ \ \ \ \ \ \ \ \ }+\left( \frac{%
\mathcal{M}_{n}}{\sqrt{n}}\right) ^{q}\left\{ \mathcal{B}_{n}(p)+\mathcal{C}%
_{n}(p)\right\} +\mathcal{D}_{n}(p)  \notag
\end{eqnarray}%
where $\mathcal{D}_{n}(p)$\ is defined as follows by case:\medskip \newline
$a.$ (\textbf{$\mathcal{L}_{rs}$-boundedness}) In general for $r$ $>$ $1$,%
\begin{equation}
\mathcal{D}_{n}(p)=2^{(r-1)/r}p^{q/(sr)}\max_{i}\left\Vert \bar{x}%
_{i,n}\right\Vert _{rs}^{q}\left( \frac{\ln (p)}{\ln \mathbb{\bar{P}}_{%
\mathcal{M}_{n}}^{-1}}\right) ^{(r-1)/r}.  \label{Dn_gen}
\end{equation}%
If additionally $\max_{i}\mathbb{E}|\bar{x}_{i,n}|^{rs}$ $=$ $O(n^{-a})$, $a$
$>$ $0$, then $\mathcal{D}_{n}(p)$ $\lesssim $ $%
2^{(r-1)/r}p^{q/(sr)}n^{-aq/(rs)}$ $\times $ $[\ln (p)/\ln \mathbb{\bar{P}}_{%
\mathcal{M}_{n}}^{-1}]^{(r-1)/r}$.$\medskip $\newline
$b.$ (\textbf{$\mathcal{L}_{rs}$-boundedness}) In general for $r$ $>$ $1$ 
\begin{equation}
\mathcal{D}_{n}(p)=2^{(r-1)/r}p^{q/(sr)}\max_{i}\left\Vert \bar{x}%
_{i,n}\right\Vert _{rs}^{q}\left( \frac{\ln (p)}{\ln \left( \mathcal{M}%
_{n}^{q}/\max_{i}\mathbb{E}\left\vert \bar{x}_{i,n}\right\vert ^{q}\right) }%
\right) ^{(r-1)/r}.  \label{Dn_Lq}
\end{equation}%
If additionally $\max_{i}\mathbb{E}|\bar{x}_{i,n}|^{rs}$ $=$ $O(n^{-a})$, $a$
$>$ $0$, then $\mathcal{D}_{n}(p)$ $\lesssim $ $%
2^{(r-1)/r}p^{q/(sr)}n^{-aq/(rs)}[\ln (p)/$ $\ln (\mathcal{M}%
_{n}^{q}n^{a})]^{(r-1)/r}$.$\medskip $\newline
$c$. (\textbf{sub-exponential}) If $\max_{i}\mathbb{P}(|\bar{x}_{i,n}|$ $%
\geq $ $c)$ $\leq $ $a\exp \{-bn^{\gamma }c^{\gamma }\}$ for some finite
constants $(a,b)$ $>$ $0$, $\gamma $ $>$ $q$, and all $c$ $>$ $0$, then for $%
p$ $>$ $e$ and $r$ $>$ $1$%
\begin{equation}
\mathcal{D}_{n}(p):=\sqrt{\frac{8}{n^{\gamma /q}a^{1/q}}\frac{\ln (p)}{\exp
\left\{ \exp \left\{ \mathcal{M}_{n}^{q}\right\} \right\} }}.
\label{Dn_sube}
\end{equation}
\end{proposition}

\begin{remark}
\normalfont Error $\mathcal{B}_{n}(p)$ arises from a Gaussian comparison
with and without truncation. $\mathcal{C}_{n}(p)$ captures the total
Gaussian approximation error from truncation and blocking. $\mathcal{D}%
_{n}(p)$ represents the truncation tail approximation error caused by
replacing $\mathbb{E}\max_{i}|\bar{x}_{i,n}|^{q}$ with $\mathbb{E}\max_{i}|%
\bar{x}_{i,n}^{(\mathcal{M})}|^{q}$.
\end{remark}

\begin{remark}
\normalfont$\mathbb{E}|\bar{x}_{i,n}|^{rs}$ $=$ $O(n^{-a})$ holds for many
stochastic processes, including $\alpha ,\beta ,\mathcal{C},\phi ,\tau
,\varphi $-mixing, mixingale, near epoch dependent, and physical dependent
processes. See \cite{McLeish1975} , \cite{Hansen1991}, \cite{Wu2005} and 
\cite{Hill2024_maxlln} to name a few.
\end{remark}

\begin{remark}
\label{rm:Mn}\normalfont The proof exploits the decomposition 
\begin{eqnarray}
\mathbb{E}\max_{i}\left\vert \bar{x}_{i,n}\right\vert ^{q} &=&\mathbb{E}%
\max_{i}\left\vert \bar{x}_{i,n}\right\vert ^{q}\mathcal{I}_{\max_{i}|\bar{x}%
_{i,n}|\leq \mathcal{M}_{n}}+\mathbb{E}\max_{i}\left\vert \bar{x}%
_{i,n}\right\vert ^{q}\mathcal{I}_{\max_{i}|\bar{x}_{i,n}|>\mathcal{M}_{n}}
\label{J1J2} \\
\text{ \ \ } &&\text{\ }=:\mathfrak{I}_{n,1}+\mathfrak{I}_{n,2}.  \notag
\end{eqnarray}%
Thus $\mathcal{M}_{n}$ imparts two opposing forces. First, \emph{large} $%
\mathcal{M}_{n}$ yields a better truncation approximation $\mathbb{E}%
\max_{i}\left\vert \bar{x}_{i,n}\right\vert ^{q}$ $\times $ $\mathcal{I}%
_{\max_{i}|\bar{x}_{i,n}|\leq \mathcal{M}_{n}}$ $\approx $ $\mathbb{E}%
\max_{i}\left\vert \bar{x}_{i,n}\right\vert ^{q}$, hence smaller $\{\mathcal{%
B}_{n}(p),\mathcal{D}_{n}(p)\}$. Second, under dependence we cannot use a
symmetrization argument, thus we use Gaussian approximations $\{\rho
_{n},\rho _{n}^{\ast (\mathcal{M})}\}$. Hence $\mathbb{E}\max_{i}\left\vert 
\bar{x}_{i,n}\right\vert ^{q}\mathcal{I}_{\max_{i}|\bar{x}_{i,n}|\leq 
\mathcal{M}_{n}}$ incurs a cost via scale $(\mathcal{M}_{n}/\sqrt{n})^{q}$:%
\begin{equation}
\mathbb{E}\max_{i}\left\vert \bar{x}_{i,n}\right\vert ^{q}\mathcal{I}%
_{\max_{i}|\bar{x}_{i,n}|\leq \mathcal{M}_{n}}\leq \mathbb{E}%
\max_{i}\left\vert \frac{1}{n}\sum_{l=1}^{\mathcal{N}_{n}}\varepsilon _{l}%
\mathcal{S}_{n,l}^{(\mathcal{M})}(i)\right\vert ^{q}+\left( \frac{\mathcal{M}%
_{n}}{\sqrt{n}}\right) ^{q}\left\{ \rho _{n}\vee \rho _{n}^{\ast (\mathcal{M}%
)}\right\} .  \notag
\end{equation}%
A \emph{smaller} $\mathcal{M}_{n}$ is then optimal. Cf. (\ref{EEx}) in the
proof of Proposition \ref{prop:main}, and step 1 in the proof of Proposition %
\ref{prop:main_unbound}.
\end{remark}

Arguments leading to (\ref{dM_U}) combined with the proof of Proposition \ref%
{prop:main} yield $\mathfrak{I}_{n,1}$ $\leq $ $\mathcal{A}_{n}(p)$ $+$ $%
\mathcal{B}_{n}(p)$ $\vee $ $\mathcal{C}_{n}(p)$. In order to bound $%
\mathfrak{I}_{n,2}$ in (\ref{J1J2}), under $\mathcal{L}_{rs}$-boundedness by
H\"{o}lder and Lyapunov inequalities, for some $s$ $>$ $q$ and $r$ $>$ $1$,%
\begin{eqnarray}
\mathfrak{I}_{n,2} &\leq &\mathbb{E}\max_{i}\left\vert \bar{x}%
_{i,n}\right\vert ^{q}\mathcal{I}_{\max_{i}|\bar{x}_{i,n}|>\mathcal{M}_{n}} 
\notag \\
&\leq &\left( \mathbb{E}\max_{i}\left\vert \bar{x}_{i,n}\right\vert
^{rq}\right) ^{1/r}\times \mathbb{P}\left( \max_{i}\left\vert \bar{x}%
_{i,n}\right\vert >\mathcal{M}_{n}\right) ^{(r-1)/r}  \notag \\
&\leq &p^{q/(sr)}\max_{i}\left( \left( \mathbb{E}\left\vert \bar{x}%
_{i,n}\right\vert ^{rs}\right) ^{1/s}\right) ^{q/r}\times \mathbb{P}\left(
\max_{i}\left\vert \bar{x}_{i,n}\right\vert >\mathcal{M}_{n}\right)
^{(r-1)/r}  \notag \\
&=&p^{q/(sr)}\max_{i}\left\Vert \bar{x}_{i,n}\right\Vert _{rs}^{q}\times 
\mathbb{P}\left( \max_{i}\left\vert \bar{x}_{i,n}\right\vert >\mathcal{M}%
_{n}\right) ^{(r-1)/r}.  \label{L4}
\end{eqnarray}%
Proposition \ref{prop:main_unbound} uses the following concentration
inequality to deal with $\mathbb{P}(\max_{i}|\bar{x}_{i,n}|$ $>$ $\mathcal{M}%
_{n})$ in (\ref{L4}). Write as before $\mathbb{\bar{P}}_{\mathcal{M}}$ $:=$ $%
\max_{i}\mathbb{P}(|\bar{x}_{i,n}|$ $>$ $\mathcal{M})$. The following result
appears new, and is akin to a tail measure version of \cite{Nemirovski2000}%
's inequality, in most cases without need to define dependence.

\begin{lemma}
\label{lm:concen}Let $\{x_{t}\}$, $x_{t}$ $=$ $[x_{i,t}]_{i=1}^{p}$, be
random variables on a common probability space.$\medskip $\newline
$a.$ If $p$ $>$ $e$ and $p=o(\mathbb{\bar{P}}_{\mathcal{M}_{n}}^{-1})$\ then 
\begin{equation}
\mathbb{P}\left( \max_{i}\left\vert \bar{x}_{i,n}\right\vert >\mathcal{M}%
_{n}\right) \leq 2\frac{\ln (p)}{\ln \mathbb{\bar{P}}_{\mathcal{M}_{n}}^{-1}}%
\rightarrow 0.  \label{PlnP}
\end{equation}%
$b.$ If $x_{i,t}$ are $\mathcal{L}_{q}$-bounded, $q$ $\geq $ $1$, $\max_{i}||%
\bar{x}_{i,n}||_{q}$ $=$ $o(\mathcal{M}_{n})$, and $p$ $=$ $o(\mathcal{M}%
_{n}^{q}/\max_{i}\mathbb{E}\left\vert \bar{x}_{i,n}\right\vert ^{q})$ then 
\begin{equation}
\mathbb{P}\left( \max_{i}\left\vert \bar{x}_{i,n}\right\vert >\mathcal{M}%
_{n}\right) \leq 2\frac{\ln (p)}{\ln \left( \mathcal{M}_{n}^{q}/\max_{i}%
\mathbb{E}\left\vert \bar{x}_{i,n}\right\vert ^{q}\right) }\rightarrow 0.
\label{Plq}
\end{equation}%
$c.$ If $\max_{i}\mathbb{P}(|\bar{x}_{i,n}|$ $\geq $ $c)$ $\leq $ $a\exp
\{-bn^{\gamma }c^{\gamma }\}$ for some $(a,b,\gamma )$ $>$ $0$ and all $c$ $%
> $ $0$, then for any $\phi $ $\in $ $(0,\gamma )$, $p$ $>$ $e$ and $\ln (p)$
$=$ $o(n^{\phi }\mathcal{M}_{n}^{\phi })$,%
\begin{equation*}
\mathbb{P}\left( \max_{i}\left\vert \bar{x}_{i,n}\right\vert \geq \mathcal{M}%
_{n}\right) \leq \frac{\ln (p)}{n^{\phi }\mathcal{M}_{n}^{\phi }\ln (\ln p)}%
+a\frac{\left( \ln (p)\right) ^{n^{\phi }\mathcal{M}_{n}^{\phi }}}{n^{\phi }%
\mathcal{M}_{n}^{\phi }\exp \left\{ bn^{\gamma }\mathcal{M}_{n}^{\gamma
}\right\} \ln (\ln p)}\rightarrow 0.
\end{equation*}
\end{lemma}

\begin{remark}
\normalfont($a$) optimizes a log-exp bound without using any tail
information, while ($b$) uses the same logic with Markov's inequality and $%
\mathcal{L}_{q}$-boundedness.
\end{remark}

\begin{remark}
\normalfont\label{rm:Bern}Sub-exponential tails for the sample mean $\max_{i}%
\mathbb{P}(|\bar{x}_{i,n}|$ $\geq $ $c)$ $\leq $ $a\exp \{-bn^{\gamma
}c^{\gamma }\}$ in ($c$) is essentially a Bernstein or Fuk-Naegev-type
inequality, and thus must be verified by case. It holds when $x_{i,t}$ has
sub-exponential tails and, e.g., satisfies physical dependence %
\citep[Theorem 2($ii$)]{Wu2005}, geometric $\tau $-mixing 
\citep[Theorem
1]{Merlevede_etal2011}, or $\alpha $-mixing or a mixingale property %
\citep{Hill2024_maxlln,Hill2025_mixg}. See \cite{Rio2017}, and see, e.g., %
\citet{Hill2024_maxlln,Hill2025_mixg} for a comprehensive review.
\end{remark}

\begin{remark}
\normalfont We still need the $\mathcal{L}_{q}$ norm to be bounded $%
\max_{i}||\bar{x}_{i,n}||_{q}$ $=$ $o(\mathcal{M}_{n})$ under $(b)$. This is
mild or irrelevant considering under common regularity conditions $\max_{i}||%
\bar{x}_{i,n}||_{q}$ $=$ $O(1/\sqrt{n})$ automatically for possibly
non-stationary mixingales and thus $\alpha $-mixing sequences %
\citep{Hansen1991,Hansen1992}, as well as physical dependent and other
mixing sequences \citep[see, e.g.,][]{Hill2024_maxlln}. Thus $\max_{i}||\bar{%
x}_{i,n}||_{q}$ $=$ $o(\mathcal{M}_{n})$ allows for trend: e.g. $x_{i,t}$ $=$
$f_{n}$ $+$ $z_{i,t}$ where $z_{i,t}$ is stationary, uniformly $\mathcal{L}%
_{q}$-bounded, and $\{f_{n}\}$ is a sequence of constants, $f_{n}$ $=$ $o(%
\mathcal{M}_{n})$. Then $\max_{i,t}||x_{i,t}/\mathcal{M}_{n}||_{q}$ $=$ $%
o(1) $ when $\mathcal{M}_{n}$ $\rightarrow $ $\infty $,\ and therefore $%
\max_{i}||\bar{x}_{i,n}||_{q}$ $=$ $o(\mathcal{M}_{n})$.
\end{remark}

\section{Examples\label{sec:Exs}}

We now study mixing and physical dependence settings in which the Gaussian
approximation Kolmogorov distances $\{\rho _{n},\rho _{n}^{\ast },\rho
_{n}^{\ast (\mathcal{M})}\}$ $\rightarrow $ $0$ yielding the error $\mathcal{%
C}_{n}(p)$ $\rightarrow $ $0$ in (\ref{BnCn}) once $p$ is bounded. We begin
with a bounded process. Throughout the multiplier $\varepsilon _{l}$ is iid,
bounded, and $\mathbb{E}\varepsilon _{l}^{2}$ $=$ $1$.

\subsection{Bounded, geometric mixing\label{sec:Ex1}}

Assume $\{x_{t}\}$ are zero mean $[-\mathcal{M}_{n},\mathcal{M}_{n}]^{p}$%
-valued random variables. Suppose $x_{i,t}$ $=$ $g_{i,t}(\epsilon _{i,t},$ $%
\epsilon _{i,t-1},\ldots )$ for some measurable functions $g_{i,t}(\cdot )$
and iid sequences $\{\epsilon _{i,t}\}$.

Define a functional $\psi _{\zeta }(x)$ $:=$ $\exp \{\left\vert x\right\vert
^{\zeta }\}$ $-$ $1$, $\zeta $ $>$ $0$, and an exponential Orlicz norm $%
||x||_{\psi _{\zeta }}$ $:=$ $\inf \{\lambda $ $>$ $0$ $:$ $\mathbb{E}[\psi
_{\zeta }\left( x/\lambda \right) ]$ $\leq $ $1\}$. Now let $%
||x_{i,t}||_{\psi _{\zeta }}$ $\leq $ $\chi _{n}$ for some $\zeta $ $\geq $ $%
1$, all $t\in \{1,...,n\}$ and $i\in \{1,...,p\}$, and some sequence $\{\chi
_{n}\}$, $\inf_{n\in \mathbb{N}}\chi _{n}$ $\geq $ $1$. Thus $\mathbb{E}\exp
\{\left\vert x_{i,t}\right\vert ^{\zeta }/\chi _{n}^{\zeta }\}$ $\leq $ $2$
implying heterogeneous sub-exponential tails which holds automatically under
boundedness. Indeed $\max_{i,t}\mathbb{E}\left\vert x_{i,t}\right\vert ^{q}$ 
$\leq $ $\mathcal{M}_{n}^{q}$. However, ignoring the upper bound $\mathcal{M}%
_{n}$ and for general use below, use the heterogeneous sub-exponential
implication $\mathbb{P}(|x_{i,t}|$ $>$ $u)$ $\leq $ $2\exp \{-u^{\zeta }\chi
_{n}^{-\zeta }\}$ with $\mathbb{E}|x_{i,t}|^{q}$ $=$ $\int_{0}^{\infty }%
\mathbb{P}(|x_{i,t}|$ $>$ $u^{1/q})du$ to deduce 
\begin{equation}
\max_{i,t}\mathbb{E}\left\vert x_{i,t}\right\vert ^{q}\leq 4\left[ \frac{q}{%
\zeta }\right] _{+}!\chi _{n}^{q}\lesssim q^{q}\chi _{n}^{q}\text{ \ }%
\forall q\geq 1\text{.}  \label{Exq}
\end{equation}

Now define $\mathcal{F}_{i,n,s}^{t}$ $:=$ $\sigma (x_{i,\tau }$ $:$ $1$ $%
\leq $ $s$ $\leq $ $\tau $ $\leq $ $t$ $\leq $ $n)$ and $\alpha $-mixing
coefficients,%
\begin{equation*}
\alpha _{i,j,n}(m):=\sup_{1\leq t\leq n}\sup_{\mathcal{C}\in \mathcal{F}%
_{i,n,-\infty }^{t},\mathcal{D}\in \mathcal{F}_{j,n,t+m}^{\infty
}}\left\vert \mathbb{P}\left( \mathcal{C}\cap \mathcal{D}\right) -\mathbb{P}(%
\mathcal{C})\mathbb{P}(\mathcal{D})\right\vert \text{ and }\alpha
_{i,n}(m):=\alpha _{i,i,n}(m).
\end{equation*}%
Assume geometric decay $\lim \sup_{n\rightarrow \infty }\max_{i,j}\alpha
_{i,j,n}(m)$ $\leq $ $a\omega ^{m}$ for some $a$ $\geq $ $1$ and $\omega $ $%
\in $ $[1/e,1)$. Lastly, assume non-degeneracy $\mathbb{E}\mathcal{X}%
_{n}^{2}(i),\mathbb{E}\mathcal{X}_{n}^{\ast 2}(i)$ $>$ $K$ for each $i,n$
and some $K$ $>$ $0$. The latter variances are finite under mixing: $\max_{i}%
\mathbb{E}\mathcal{X}_{n}^{2}(i)$ $\lesssim $ $\max_{i,t}%
\{||x_{i,t}||_{r}^{2}\}$ $<$ $\infty $ for $r$ $>$ $q$ follows from %
\citet[Corollary]{Davydov1968} and geometric decay. The same argument yields 
$\max_{i,l}\mathbb{E}(1/\sqrt{b_{n}}%
\sum_{t=(l-1)b_{n}+1}^{lb_{n}}x_{i,t})^{2}$ $\lesssim $ $\max_{i,t}%
\{||x_{i,t}||_{r}^{2}\}$, hence $\max_{i}\mathbb{E}\mathcal{X}_{n}^{\ast
2}(i)$ $<$ $\infty $\ by mutual independence and $\mathbb{E}\varepsilon
_{l}^{2}$ $=$ $1$.

The use of telescoping blocks and therefore block-wise dependent $\{\eta
_{t}\}$ lowers the mixing decay rate, hampering high dimensional theory %
\citep[see][Appendix, Proposition 3]{ChangJiangShao2023}. We can yield
better bounds on $\{\rho _{n},\rho _{n}^{\ast }\}$, however, by assuming $%
x_{i,t}$ is a measurable function of iid innovations and using a new link
between mixing and physical dependence properties %
\citep{Hill2024_maxlln,Hill2025_mixg}.

\begin{lemma}
\label{lm:ex1}Assume block size $b_{n}$ and tail heterogeneity $\chi _{n}$
satisfy $\chi _{n}$ $=$ $o(b_{n})$. Let $\ln (p)$ $<$ $n^{1/4}/b_{n}$ $%
\rightarrow $ $\infty $. Then 
\begin{eqnarray}
&&\rho _{n}\lesssim \frac{\chi _{n}\left( \ln (p)\right) ^{7/6}}{n^{1/9}}
\label{rex1} \\
&&\rho _{n}^{\ast }\lesssim \left\{ \left( \frac{\chi _{n}}{b_{n}}\right)
^{1/3}\times \left( 1+\ln (p)+\left\vert \ln \left( \frac{\chi _{n}}{b_{n}}%
\right) \right\vert \right) ^{2/3}\right\} \vee \frac{\chi _{n}b_{n}^{\beta
^{\ast }}\left( \ln (p)\right) ^{\beta ^{\ast }+2}}{n^{\beta ^{\ast
}/(4+\beta ^{\ast })}}  \label{rex1_til}
\end{eqnarray}%
where $\beta ^{\ast }$ $:=$ $\sqrt{4\ln n/[\ln b_{n}+\ln \ln (p)]}$ $-$ $4$ $%
>$ $0$ $\forall n$. We therefore achieve a \cite{Nemirovski2000}-like moment
bound under geometric mixing,%
\begin{eqnarray*}
\mathbb{E}\max_{i}\left\vert \bar{x}_{i,n}\right\vert ^{q} &\lesssim &\left( 
\frac{2\ln \left( 2p\right) }{\mathcal{N}_{n}}\right) ^{q/2}\sqrt{\mathbb{E}%
\max_{i,l}\left\vert \frac{1}{b_{n}}\mathcal{S}_{n,l}(i)\right\vert ^{q}%
\mathbb{E}\left( \max_{i}\frac{1}{\mathcal{N}_{n}}\sum_{l=1}^{\mathcal{N}%
_{n}}\left\vert \mathcal{S}_{n,l}(i)\right\vert \right) ^{q}}+o(1) \\
&=&\mathfrak{N}_{n}+o(1),
\end{eqnarray*}%
with $\mathfrak{N}_{n}$\ implicit, provided for some $r$ $>$ $1$ 
\begin{eqnarray*}
\ln (p) &=&o\left( \frac{n^{1/4}}{b_{n}}\wedge \frac{n^{3q/4}}{%
\max_{i,t}\left\Vert x_{i,t}\right\Vert _{2r}^{r}\mathcal{M}_{n}^{3q/2-(r-1)}%
}\wedge \frac{n^{3q/7+2/21}}{\chi _{n}^{6/7}\mathcal{M}_{n}^{6q/7}}\wedge
\left( \frac{b_{n}}{\chi _{n}}\right) ^{1/2}\frac{n^{3q/4}}{\mathcal{M}%
_{n}^{3q/2}}\right. \\
&&\text{ \ \ \ \ \ \ \ }\left. \wedge \frac{n^{\beta ^{\ast }/[(4+\beta
^{\ast })(\beta ^{\ast }+2)]}}{\chi _{n}^{1/(\beta ^{\ast }+2)}b_{n}^{\beta
^{\ast }/(\beta ^{\ast }+2)}}\wedge n^{\gamma /q}\exp \left\{ \exp \left\{ 
\mathcal{M}_{n}^{q}\right\} \right\} \right) .
\end{eqnarray*}
\end{lemma}

\begin{remark}
\normalfont The upper bound component $(\chi _{n}/b_{n})^{1/3}$\ in (\ref%
{rex1_til}) arises from a Gaussian-to-Gaussian comparison, thus we enforce $%
\chi _{n}/b_{n}$ $\rightarrow $ $0$, an implicit restriction on
heterogeneity. A larger/faster block size $b_{n}$ $\rightarrow $ $\infty $
naturally improves the approximation between $\mathbb{E}\mathcal{X}_{n}(i)%
\mathcal{X}_{n}(k)$ and its blocked version $\mathbb{E}\mathcal{X}_{n}^{\ast
}(i)\mathcal{X}_{n}^{\ast }(j)$, and therefore between the limiting
max-Gaussian processes for $\{\max_{i}|\mathcal{X}_{n}(i)|,\max_{i}|\mathcal{%
X}_{n}^{\ast }(i)|\}$ \citep[Theorem
2]{Chernozhukov_etal2015}. Otherwise, $b_{n}$ adversely affects $\rho
_{n}^{\ast }$ $\rightarrow $ $0$, and slows down high dimension divergence $%
p $ $\rightarrow $ $\infty $, a penalty for having dependent data. The final
term $n^{\gamma /q}\exp \left\{ \exp \left\{ \mathcal{M}_{n}^{q}\right\}
\right\} $ never dominates when $\mathcal{M}_{n}$ $\rightarrow $ $\infty $.
Hence as long as $r$ $<$ $1$ $+$ $3q/2$ (which can always be assumed),
larger/faster tail heterogeneity or upper bound $(\chi _{n},\mathcal{M}_{n})$
$\rightarrow $ $\infty $ adversely affects $\ln (p)$ $\rightarrow $ $\infty $%
, penalties for deviating from a homogeneous setting.
\end{remark}

\begin{remark}
\normalfont Boundedness is not exploited when handling $\{\rho _{n},\rho
_{n}^{\ast }\}$; cf. Remark \ref{rm:Kol_M}.
\end{remark}

\begin{remark}
\normalfont Suppose block size $b_{n}$ $\simeq $ $n^{b}$, heterogeneity $%
\chi _{n}\simeq $ $n^{\chi }$, and bound $\mathcal{M}_{n}\simeq $ $n^{m}$
for finite $(b,m,\chi )$ $>$ $0$, with $\chi $ $\in $ $[0,(1/9)\wedge b)$
and $b$ $<$ $1/4$. Let $\max_{i,t}||x_{i,t}||_{2r}$ $=$ $O(1)$. Then $\ln
(p) $ $=$ $o(n^{a})$ and $\beta ^{\ast }$ $\simeq $ $2/\sqrt{b}$ $-$ $4$ with%
\begin{eqnarray*}
a &=&\left( \frac{1}{4}-b\right) \wedge \left( \frac{3q}{4}-m\left( \frac{3q%
}{2}-(r-1)\right) \right) \wedge \frac{1}{7}\left( 3q+\frac{2}{3}-6\chi
-6qm\right) \\
&&\text{ \ \ \ \ \ \ \ \ }\wedge \frac{1}{2}\left( \frac{3q}{2}+b-\chi
-3qm\right) \wedge \frac{\sqrt{b}}{2-2\sqrt{b}}\left( 1-2\sqrt{b}-\chi -2%
\sqrt{b}\left[ 1-2\sqrt{b}\right] \right) .
\end{eqnarray*}%
For example, if $q$ $=$ $1$ and $r$ $=$ $1.5$, with bound growth $m$ $=$ $%
1/2 $, homogeneity $\chi $ $=$ $0$, and block size $b$ $=$ $1/8$\ then $\ln
(p)$ $=$ $o(n^{.02345})$.
\end{remark}

\subsection{Unbounded, geometric mixing, sub-exponential\label{sec:Ex2}}

Reconsider the Example \ref{sec:Ex1} setting, except assume $x_{t}$ has
support $\mathbb{R}^{p}$. Let $\max_{i,t}\mathbb{E}|x_{i,t}|^{2r}$ $<$ $%
\infty $ for each $n$ and some $r$ $>$ $1$, and assume $\chi _{n}$ $=$ $%
o(b_{n})$ as above. Under sub-exponentiality $\max_{i}\mathbb{P}(|x_{i,t}|$ $%
>$ $u)$ $\leq $ $2\exp \{-u^{\gamma }\chi _{n}^{-\gamma }\}$ we can make $r$
arbitrarily large, cf. (\ref{Exq}), and by Proposition \ref%
{prop:main_unbound}.c for $p$ $>$ $e$,%
\begin{eqnarray*}
\mathbb{E}\max_{i}\left\vert \bar{x}_{i,n}\right\vert ^{q} &\lesssim &\left( 
\frac{2\ln \left( 2p\right) }{\mathcal{N}_{n}}\right) ^{q/2}\sqrt{\mathbb{E}%
\max_{i,l}\left\vert \frac{1}{b_{n}}\mathcal{S}_{n,l}^{(\mathcal{M}%
)}(i)\right\vert ^{q}\mathbb{E}\left( \max_{i}\frac{1}{\mathcal{N}_{n}}%
\sum_{l=1}^{\mathcal{N}_{n}}\left\vert \mathcal{S}_{n,l}^{(\mathcal{M}%
)}(i)\right\vert \right) ^{q}} \\
&&+\frac{\mathcal{M}_{n}^{q}}{n^{q/2}}\left\{ \left( \mathcal{M}%
_{n}^{-2(r-1)}\max_{i,t}\left\Vert x_{i,t}\right\Vert _{2r}^{2r}\right)
^{1/2}\ln (p)\right\} ^{2/3} \\
&&+\left( \frac{\mathcal{M}_{n}}{\sqrt{n}}\right) ^{q}\left\{ \rho _{n}+\rho
_{n}^{\ast (\mathcal{M})}\right\} +\sqrt{\frac{8}{n^{\gamma /q}a^{1/q}}}%
\frac{\sqrt{\ln (p)}}{\exp \left\{ \exp \left\{ \mathcal{M}_{n}^{q}\right\}
/2\right\} } \\
&& \\
&=&\mathcal{A}_{n}(p)+\mathcal{B}_{n}(p)+\mathcal{C}_{n}(p)+\mathcal{D}%
_{n}(p).
\end{eqnarray*}%
Analyzing $\mathcal{C}_{n}(p)$ $=$ $(\mathcal{M}_{n}/\sqrt{n})^{q}\{\rho _{n}
$ $+$ $\rho _{n}^{\ast (\mathcal{M})}\}$ proceeds exactly as in Lemma \ref%
{lm:ex1} since the latter does not exploit boundedness, thus (\ref{rex1})
and (\ref{rex1_til}) apply. Furthermore, after simplifying terms it can be
shown all $(\mathcal{B}_{n}(p),\mathcal{C}_{n}(p),\mathcal{D}_{n}(p))$ $%
\rightarrow $ $0$ sufficiently when%
\begin{eqnarray}
\ln (p) &=&o\left( \frac{n^{3q/4}}{\max_{i,t}\left\Vert x_{i,t}\right\Vert
_{2r}^{r}\mathcal{M}_{n}^{3q/2-(r-1)}}\wedge \frac{n^{3q/7+2/21}}{\chi
_{n}^{6/7}\mathcal{M}_{n}^{6q/7}}\wedge \left( \frac{b_{n}}{\chi _{n}}%
\right) ^{1/2}\frac{n^{3q/4}}{\mathcal{M}_{n}^{3q/2}}\right.   \label{lnp0}
\\
&&\text{ \ \ \ \ \ \ \ }\left. \wedge \frac{n^{[q/2+\beta ^{\ast }/(4+\beta
^{\ast })]/(\beta ^{\ast }+2)}}{\chi _{n}^{1/(\beta ^{\ast }+2)}b_{n}^{\beta
^{\ast }/(\beta ^{\ast }+2)}\mathcal{M}_{n}^{q/(\beta ^{\ast }+2)}}\wedge
n^{\gamma /q}\exp \left\{ \exp \left\{ \mathcal{M}_{n}^{q}\right\} \right\}
\right) .  \notag
\end{eqnarray}

Finally, recall $\mathcal{M}_{n}$ satisfies (\ref{Mn}), while we need $%
\mathcal{M}_{n}^{q}n^{-q/2}\{(\mathcal{M}_{n}^{-2(r-1)}\max_{i,t}\left\Vert
x_{i,t}\right\Vert _{2r}^{2r})^{1/2}\}^{2/3}$ $\rightarrow $ $0$ in $%
\mathcal{B}_{n}(p)$. Bound heterogeneity by assuming $%
\max_{i,t}||x_{i,t}||_{2r}^{r}$ $=$ $O(1)$. Then we need $n^{1/(r-1)}$ $<$ $%
\mathcal{M}_{n}$ $<$ $n^{q/\{2\left[ q-2(r-1)/3\right] \}}$, hence $r$ $>$ $%
(4$ $+$ $9q)/(4$ $+$ $3q)$. Thus if $q$ $=$ $2$ then $r$ $>$ $2.2$, hence $%
x_{i,t}$ must be $\mathcal{L}_{4.4+\delta }$-bounded for some $\delta $ $>$ $%
0$.

\subsubsection{ Example: sub-Gaussian}

In a sub-Gaussian setting $\gamma $ $=$ $2$ thus $||\max_{i,t}\left\vert
x_{i,t}\right\vert ||_{2}$ $=$ $O(\chi _{n}^{1/2}\ln (pn)^{1/4})$. For
example, if $|x_{i,t}|^{rq}$ has sub-exponential tails then $\mathbb{P}%
(|x_{i,t}|$ $>$ $u)$ $\leq $ $2\exp \{-u^{rq\gamma }\chi _{n}^{-\gamma }\}$, 
$\gamma $ $\geq $ $1$. Thus if $rq$ $\geq $ $2/\gamma $ then $x_{i,t}$ is
sub-Gaussian; e.g. $q$ $=$ $2$, \ $r$ $=$ $1.5$ and $\gamma $ $\geq $ $2/3$.

Assume bounded heterogeneity $\chi _{n}$ $=$ $O(1)$ for simplicity, set
truncation level $\mathcal{M}_{n}$ $=$ $n^{m}$, $m$ $>$ $0$, block size $%
b_{n}$ $\simeq $ $n^{b}$, $b$ $\in $ $(0,1)$, and set $q$ $=$ $1$ and $r$ $%
=1.5$. Thus $\beta ^{\ast }$ $\simeq $ $2/\sqrt{b}$ $-$ $4$, and $\ln (p)$ $%
= $ $o(n^{a})$ where 
\begin{equation*}
a=\left( \frac{3}{4}-m\right) \wedge \left( \frac{3-6m}{7}+\frac{2}{21}%
\right) \wedge \left( \frac{b-3m}{2}+\frac{3}{4}\right) \wedge \left( \frac{1%
}{\beta ^{\ast }+2}\left[ \frac{1}{2}+\frac{\beta ^{\ast }}{4+\beta ^{\ast }}%
\right] -\frac{b\beta ^{\ast }+m}{\beta ^{\ast }+2}\right) .
\end{equation*}%
We therefore need $b$ $<$ $1/4$, and $m$ $<$ $\{b/3$ $+$ $1/2\}\wedge \{1/2$ 
$+$ $(2/\sqrt{b}$ $-$ $4)(.5\sqrt{b}$ $-$ $b)\}$.

\subsection{Unbounded, physical dependent, $\mathcal{L}_{q}$-bounded\label%
{sec:Ex3}}

Now consider physical dependence as in \cite{Wu2005} and \cite{WuMin2005}.
We work in the high dimensional Gaussian approximation setting of \cite%
{ZhangCheng2018} allowing for non-stationarity, but results in, e.g., \cite%
{ChangChenWu2024} also work.

Let $\{\epsilon _{i,t}^{\prime }\}$ be an independent copy of $\{\epsilon
_{i,t}\}$, and define the coupled process $x_{i,n,t}^{\prime }(m)$ $:=$ $%
g_{i,n,t}(\epsilon _{i,t},\ldots ,\epsilon _{i,t-m+1},\epsilon
_{i,t-m}^{\prime },\epsilon _{i,t-m-1},\ldots )$, $m$ $=$ $0,1,2,...$ Define
the $\mathcal{L}_{q}$-\textit{physical dependence} measure $\theta
_{i,n}^{(q)}(m)$ and its partial accumulation 
\begin{equation*}
\theta _{i,n}^{(q)}(m):=\max_{t}\left\Vert x_{i,t}-x_{i,t}^{\prime
}(m)\right\Vert _{q}\text{ \ and \ }\Theta
_{i,n}^{(q)}(m):=\sum_{l=m}^{\infty }\theta _{i,n}^{(q)}(l).
\end{equation*}%
We say $x_{t}$ is (uniformly) $\mathcal{L}_{q}$-\textit{physical dependent}
when $\max_{i}\Theta _{i,n}^{(q)}(0)$ $<$ $K$ for each $n$.

We also need the truncated process $x_{i,t}^{(\mathcal{M})}$ $=$ $x_{i,t}%
\mathcal{I}_{\left\vert x_{i,t}\right\vert \leq \mathcal{M}_{n}}$ $+$ $%
\mathcal{M}_{n}\mathcal{I}_{\left\vert x_{i,t}\right\vert >\mathcal{M}_{n}}$
to satisfy physical dependence. It is easy to show the property carries over
to affine combinations and products \citep[e.g.][]{Wu2005}, however it
generally need not carry over to \textit{general} measurable transforms,
including $\mathcal{I}_{\left\vert x_{i,t}\right\vert \leq \mathcal{M}_{n}}$.

In Assumption \ref{assume:ex_phys}, below, we merely assume $x_{i,t}^{(%
\mathcal{M})}$ is physical dependent. However, under mild conditions it
holds when $x_{i,t}$ is physical dependent. In order to see this define $%
\mathcal{I}_{i,n,t}$ $:=$ $\mathcal{I}_{\left\vert x_{i,t}\right\vert \leq 
\mathcal{M}_{n}}$, a coupled version $\mathcal{I}_{i,n,t}^{\prime }(m)$ $:=$ 
$\mathcal{I}_{|x_{i,t}^{\prime }(m)|\leq \mathcal{M}_{n}}$, and coefficients%
\begin{eqnarray*}
&&\varpi _{i,n}^{(q)}(m):=\max_{t}\left\Vert \mathcal{I}_{i,n,t}-\mathcal{I}%
_{i,n,t}^{\prime }(m)\right\Vert _{q} \\
&&\theta _{i,n}^{\mathcal{M}(q)}(m):=\max_{t}\left\Vert y_{i,t}^{(\mathcal{M}%
)}-y_{i,t}^{(\mathcal{M})\prime }(m)\right\Vert _{q}\text{ \ and \ }\Theta
_{i,n}^{\mathcal{M}(q)}(m):=\sum_{l=m}^{\infty }\theta _{i,n}^{\mathcal{M}%
(q)}(l).
\end{eqnarray*}%
By Minkowski, Cauchy-Schwartz and Lyapunov inequalities, for any $q$ $\geq $ 
$2$%
\begin{eqnarray*}
\theta _{i,n}^{\mathcal{M}(q/2)}(m) &\leq &\max_{t}\left\Vert x_{i,t}%
\mathcal{I}_{\left\vert x_{i,t}\right\vert \leq \mathcal{M}%
_{n}}-x_{i,t}^{\prime }(m)\mathcal{I}_{|x_{i,t}^{\prime }(m)|\leq \mathcal{M}%
_{n}}\right\Vert _{q/2} \\
&&\text{ \ \ \ }+\mathcal{M}_{n}\max_{t}\left\Vert \mathcal{I}_{\left\vert
x_{i,t}\right\vert >\mathcal{M}_{n}}-\mathcal{I}_{x_{i,t}^{\prime }(m)>%
\mathcal{M}_{n}}\right\Vert _{q/2} \\
&& \\
&\leq &\max_{t}\left\Vert \left( x_{i,t}-x_{i,t}^{\prime }(m)\right) 
\mathcal{I}_{\left\vert x_{i,t}\right\vert \leq \mathcal{M}_{n}}\right\Vert
_{q/2} \\
&&\text{ \ \ \ }+\max_{t}\left\Vert x_{i,t}^{\prime }(m)\mathcal{I}%
_{\left\vert x_{i,t}\right\vert \leq \mathcal{M}_{n}}-x_{i,t}^{\prime }(m)%
\mathcal{I}_{|x_{i,t}^{\prime }(m)|\leq \mathcal{M}_{n}}\right\Vert _{q/2} \\
&&\text{ \ \ \ }+\mathcal{M}_{n}\max_{t}\left\Vert \mathcal{I}_{\left\vert
x_{i,t}\right\vert >\mathcal{M}_{n}}-\mathcal{I}_{x_{i,t}^{\prime }(m)>%
\mathcal{M}_{n}}\right\Vert _{q/2} \\
&& \\
&\leq &\theta _{i,n}^{(q/2)}(m)+\left( \max_{t}\left\Vert x_{i,t}^{\prime
}(m)\right\Vert _{q}+\mathcal{M}_{n}\right) \times \varpi _{i,n}^{(q)}(m) \\
&\leq &\theta _{i,n}^{(q/2)}(m)+\left( \max_{t}\left\Vert x_{i,t}\right\Vert
_{q}+\mathcal{M}_{n}\right) \times \varpi _{i,n}^{(q)}(m).
\end{eqnarray*}%
Moreover, for any $q$ $>$ $1$ and $m$ $\geq $ $0$, 
\begin{equation*}
\varpi _{i,n}^{(q)}(m):=\max_{t}\left\Vert \mathcal{I}_{i,n,t}-\mathcal{I}%
_{i,n,t}^{\prime }(m)\right\Vert _{q}\leq \max_{t}\mathbb{P}\left(
\left\vert x_{i,t}\right\vert >\mathcal{M}_{n}\right) \leq \mathcal{M}%
_{n}^{-q}\max_{t}\left\Vert x_{i,t}\right\Vert _{q}^{q}.
\end{equation*}%
Thus%
\begin{equation}
\theta _{i,n}^{\mathcal{M}(q/2)}(m)\leq \theta _{i,n}^{(q/2)}(m)+\left( 
\frac{\max_{t}\left\Vert x_{i,t}\right\Vert _{q}}{\mathcal{M}_{n}^{q/(q+1)}}%
\right) ^{q+1}+\left( \frac{\max_{t}\left\Vert x_{i,t}\right\Vert _{q}}{%
\mathcal{M}_{n}^{1-1/q}}\right) ^{q}.  \label{th_M}
\end{equation}%
Since $q/(q$ $+$ $1)$ $>$ $1$ $-$ $1/q$ we conclude $x_{i,t}^{(\mathcal{M})}$
is $\mathcal{L}_{q/2}$-physical dependent when $x_{i,t}$ is, as long as $%
\max_{t}||x_{i,t}||_{q}$ $=$ $o(\mathcal{M}_{n}^{1-1/q})$. Thus physical
dependence extends to $x_{i,t}^{(\mathcal{M})}$ as long as trend is not too
severe, with a trade-off with the availability of higher moments $q$.

Write $l_{n}$ $:=$ $l_{n}(p,\gamma )$ $=$ $\ln (pn/\gamma _{n})\vee 1$ for
any positive sequence $\{\gamma _{n}\}$.

\begin{assumption}
\label{assume:ex_phys} $\ \ \ \medskip $\newline
$a.$ Let $\max_{i,t}\mathbb{E}x_{i,t}^{4}$ $\leq $ $K$; $\max_{t}\mathbb{E}%
h(\max_{i}|x_{i,t}|/\chi _{n})$ $\leq $ $K$ for some strictly increasing
convex function $h$ $:$ $[0,\infty )$ $\rightarrow $ $[0,\infty )$, and
sequence $\{\chi _{n}\},$ $\chi _{n}$ $\geq $ $1$.\medskip \newline
$b.$ There exist $\mathcal{E}_{n}$ $>$ $0$, $\mathcal{E}_{n}$ $\rightarrow $ 
$\infty $, and $\gamma _{n}$ $\in $ $(0,1)$ such that $K\{\chi
_{n}h^{-1}\left( n/\gamma _{n}\right) $ $\vee $ $l_{n}^{1/2}\}$ $\leq $ $%
n^{3/8}/[\mathcal{E}_{n}^{1/2}l_{n}^{5/8}]$.$\medskip $\newline
$c.$ $0$ $<$ $\min_{i}\mathbb{E}\mathcal{X}_{n}^{2}(i)$ $\leq $ $\max_{i}%
\mathbb{E}\mathcal{X}_{n}^{2}(i)$ $<$ $K$ for each $n$; $\{x_{i,t},x_{i,t}^{(%
\mathcal{M})}\}$ are uniformly physical dependent with $\max_{i}\sum_{l=0}^{%
\infty }\theta _{i,n}^{\mathcal{M}(4)}(l)$ $<$ $\infty $ and $%
\max_{i}\sum_{l=0}^{\infty }l\{\theta _{i,n}^{(3)}(l)$ $\vee $ $\theta
_{i,n}^{\mathcal{M}(3)}(m)\}$ $\leq $ $\xi _{n}$ for some sequence $\xi _{n}$
$\geq $ $1$.
\end{assumption}

\begin{remark}
\normalfont An exponential map $h(z)$ $=$ $\exp \{z\}$ implies $x_{i,t}$
requires sub-exponentail tails. The $l_{q}$ map $h(z)$ $=$ $z^{q}$, $q$ $%
\geq $ $1$, implies $h^{-1}(y)$ $=$ $y^{1/q}$, hence ($b$) reduces to%
\begin{equation*}
\left\{ \chi _{n}(n/\gamma _{n})^{1/q}\vee \left( \ln (pn/\gamma _{n})\vee
1\right) ^{1/2}\right\} \lesssim \frac{n^{3/8}}{\mathcal{E}_{n}^{1/2}\left(
\ln (pn/\gamma _{n})\vee 1\right) ^{5/8}}.
\end{equation*}
\end{remark}

\begin{remark}
\normalfont Assumption \ref{assume:ex_phys} yields Assumptions 2.1-2.3 in 
\cite{ZhangCheng2018}. The generalization $\max_{i}\sum_{l=1}^{\infty
}l\{\theta _{i,n}^{(6)}(l)$ $\vee $ $\theta _{i,n}^{\mathcal{M}(3)}(m)\}$ $%
\leq $ $\xi _{n}$ allows for trending moments that may not be adequately
captured in $\max_{t}\mathbb{E}h(\max_{i}|x_{i,t}|/\chi _{n})$ $\leq $ $K$
with $\chi _{n}$ for a given $h$. If, however, $h(\cdot )$ $=$ $\exp \{\cdot
\}$ or $|\cdot |^{q}$ then $\max_{i,t}\mathbb{E}|x_{i,t}|^{q}$ $\lesssim $ $%
\chi _{n}^{q}$ by (\ref{Exq}), hence we may take $\xi _{n}$ $\simeq $ $\chi
_{n}$.
\end{remark}

\begin{remark}
\normalfont For the first part of ($c$), by arguments in 
\citet[Theorems 1
and 2(i)]{Wu2005} and the second part of ($c$), $\max_{i}\mathbb{E}\mathcal{X%
}_{n}^{2}(i)$ $\lesssim $ $\max_{i}\sum_{l=1}^{\infty }\theta _{i,n}^{(2)}(l)
$ $\leq $ $K$. Verification of the second part of ($c$) requires either
additional dependence or heterogeneity restrictions because $f(x_{i,t})$
need not be physical dependent for measurable $f$. ($i$) Suppose $x_{i,t}$
are $\alpha $-mixing with $\max_{i.j}\alpha _{i,j,n}(m)$ $\lesssim $ $%
m^{-\lambda }$ for some $\lambda $ $>$ $2q^{2}/(q$ $-$ $1)$ and $q$ $\geq $ $%
6$. Then $\{\mathcal{I}_{\left\vert x_{i,t}\right\vert \leq \mathcal{M}%
_{n}},x_{i,t}\mathcal{I}_{\left\vert x_{i,t}\right\vert \leq \mathcal{M}%
_{n}}\}$ are $\alpha $-mixing with coefficients $\alpha _{i,n}(m)$, and
hence $\{x_{i,t},x_{i,t}^{(\mathcal{M})}\}$ are physical dependent with
coefficients $||x_{i,t}||_{q}^{1/q}\alpha _{i,n}^{1/(qr)}(m)$ $\lesssim $ $%
||x_{i,t}||_{q}^{1/q}m^{-\lambda (q-1)/q^{2}}$ where $\lambda (q$ $-$ $%
1)/q^{2}$ $>$ $2$: see \citet[Lemma 1.6]{McLeish1975} and \citet[Theorem
2.1]{Hill2025_mixg}, cf. \citet[proof of Theorem 2.10]{Hill2024_maxlln}.
Hence $\max_{i}\sum_{l=1}^{\infty }l\{\theta _{i,n}^{(3)}(l)$ $\vee $ $%
\theta _{i,n}^{\mathcal{M}(3)}(m)\}$ $\leq $ $K$, where $q$ $>$ $1$ must
hold, a typical \textquotedblleft \emph{moment-memory}\textquotedblright\
trade-off. Or: ($ii$) Let $\max_{i}\sum_{l=0}^{\infty }l\theta
_{i,n}^{(3)}(l)$ $\leq $ $\xi _{n}$ and $\max_{t}||x_{i,t}||_{q}$ $=$ $o(%
\mathcal{M}_{n}^{1-1/q})$. Then $\max_{i}\sum_{l=0}^{\infty }l\theta _{i,n}^{%
\mathcal{M}(3)}(l)$ $\leq $ $\xi _{n}$ by the arguments leading to (\ref%
{th_M}).
\end{remark}

Under $\mathcal{L}_{q}$-boundedness and Proposition \ref{prop:main_unbound}%
.b,%
\begin{eqnarray}
\text{ \ \ \ }\mathbb{E}\max_{i}\left\vert \bar{x}_{i,n}\right\vert ^{q}
&\lesssim &\left( \frac{2\ln \left( 2p\right) }{\mathcal{N}_{n}}\right)
^{q/2}\sqrt{\mathbb{E}\max_{i,l}\left\vert \frac{1}{b_{n}}\mathcal{S}%
_{n,l}^{(\mathcal{M})}(i)\right\vert ^{q}\mathbb{E}\left( \max_{i}\frac{1}{n}%
\sum_{l=1}^{\mathcal{N}_{n}}\left\vert \mathcal{S}_{n,l}^{(\mathcal{M}%
)}(i)\right\vert \right) ^{q}}  \label{Emax} \\
&&\text{ \ \ \ \ \ \ \ }+\left( \frac{\mathcal{M}_{n}}{\sqrt{n}}\right)
^{q}\left\{ \left( \mathcal{M}_{n}^{-2(r-1)}\max_{i,t}\left\Vert
x_{i,t}\right\Vert _{2r}^{2r}\right) ^{1/2}\ln (p)\right\} ^{2/3}  \notag \\
&&\text{ \ \ \ \ \ \ \ }+\left( \frac{\mathcal{M}_{n}}{\sqrt{n}}\right)
^{q}\left\{ \rho _{n}+\rho _{n}^{\ast (\mathcal{M})}\right\}  \notag \\
&&\text{ \ \ \ \ \ \ \ }+\left\Vert \max_{i,t}\left\vert x_{i,t}\right\vert
\right\Vert _{rq}^{q}\left( \frac{\ln (p)}{\ln \left( \mathcal{M}%
_{n}^{q}/\max_{i}\mathbb{E}\left\vert \bar{x}_{i,n}\right\vert ^{q}\right) }%
\right) ^{(r-1)/r}  \notag \\
&=&\mathcal{A}_{n}(p)+\mathcal{B}_{n}(p)+\mathcal{C}_{n}(p)+\mathcal{D}%
_{n}(p).  \notag
\end{eqnarray}%
It remains to bound $\{\rho _{n},\rho _{n}^{\ast (\mathcal{M})}\}$ in $%
\mathcal{C}_{n}(p)$. Recall $l_{n}$ $:=$ $l_{n}(p,\gamma )$ $=$ $\ln
(pn/\gamma _{n})\vee 1$ for some $\gamma _{n}$ $\in $ $(0,1)$, and define%
\begin{eqnarray*}
&&\Theta _{i,n}^{(q)}(m):=\sum_{l=m}^{\infty }\left\{ \theta
_{i,n}^{(2q)}(l)\vee \theta _{i,n}^{(q)}(m)\right\} \text{ \ and \ }\Xi
_{n}(m):=\max_{i}\sum_{l=m}^{\infty }l\left\{ \theta _{i,n}^{(4)}(l)\vee
\theta _{i,n}^{(2)}(m)\right\} \\
&&\check{\Theta}_{i,n}^{(q)}(m):=\sum_{l=m}^{\infty }\left\{ \theta
_{i,n}^{(2q)}(l)\vee \theta _{i,n}^{\mathcal{M}(q)}(m)\right\} \text{ \ and
\ }\check{\Xi}_{n}(m):=\max_{i}\sum_{l=m}^{\infty }l\left\{ \theta
_{i,n}^{(4)}(l)\vee \theta _{i,n}^{\mathcal{M}(2)}(m)\right\} .
\end{eqnarray*}%
Notice $\theta _{i,n}^{\mathcal{M}(\cdot )}(m)$ under truncation in the
latter $\{\check{\Theta}_{i,n}^{(q)}(m),\check{\Xi}_{n}(m)\}$. Recall $%
\mathcal{E}_{n}$ $\rightarrow $ $\infty $ in Assumption \ref{assume:ex_phys}%
.b.

\begin{lemma}
\label{lm:ex_phys}Under Assumption \ref{assume:ex_phys}, $\xi _{n}$ $=$ $%
o(b_{n})$ and $q$ $\geq $ $2$,%
\begin{eqnarray}
&&a.\text{ }\rho _{n}\lesssim \frac{\mathcal{E}_{n}^{1/2}l_{n}^{7/8}}{n^{1/8}%
}+\gamma _{n}+\left( \frac{n^{1/8}}{\mathcal{E}_{n}^{1/2}l_{n}^{3/8}}\right)
^{q/(q+1)}\left( \sum_{i=1}^{p}\Theta _{i,n}^{(q)}(\mathcal{E}_{n})\right)
^{1/(q+1)}  \label{rhonM} \\
&&\text{ \ \ \ \ \ \ \ \ \ \ \ \ \ \ \ }+\Xi _{n}^{1/3}(\mathcal{E}%
_{n})\times \left( 1\vee \ln \left( p/\Xi _{n}(\mathcal{E}_{n})\right)
\right) ^{2/3}  \notag \\
&&  \notag \\
&&b.\text{ }\rho _{n}^{\ast (\mathcal{M})}\lesssim c_{n}(p)\vee \left( \frac{%
\xi _{n}}{b_{n}}\right) ^{1/3}\ln (p)^{2/3},  \label{rhon_til}
\end{eqnarray}%
where 
\begin{eqnarray*}
&&c_{n}(p):=\frac{\mathcal{E}_{n}^{1/2}l_{n}^{7/8}}{n^{1/8}}+\gamma
_{n}+\left( \frac{n^{1/8}}{\mathcal{E}_{n}^{1/2}l_{n}^{3/8}}\right)
^{q/(q+1)}\left( \sum_{i=1}^{p}\check{\Theta}_{i,n}^{(\mathcal{M})(q/2)}(%
\mathcal{E}_{n})\right) ^{1/(q+1)} \\
&&\text{ \ \ \ \ \ \ \ \ \ \ \ \ \ \ \ \ \ }+\check{\Xi}_{n}^{(\mathcal{M}%
)1/3}(\mathcal{E}_{n})\times \left( 1\vee \ln \left( p/\check{\Xi}_{n}^{(%
\mathcal{M})}(\mathcal{E}_{n})\right) \right) ^{2/3}.
\end{eqnarray*}
\end{lemma}

\subsubsection{Example: $l_{q}$ map and hyperbolic memory}

Consider $h(x)$ $=$ $|x|^{q}$, $q$ $\geq $ $2$. Assume $\max_{i,t}\mathbb{E}%
x_{i,t}^{4}$ $\leq $ $K$, $\max_{t}\mathbb{E}\max_{i}|x_{i,t}/\chi _{n}|^{q}$
$\leq $ $K$ and $\ln (p)$ $\lesssim $ $n^{1/q}$. Assume hyperbolic weak
dependence $\sum_{i=1}^{p}\{\theta _{i,n}^{(6)}(m)$ $\vee $ $\theta _{i,n}^{%
\mathcal{M}(3)}(m)\}$ $\leq $ $\phi _{n}m^{-\lambda }$ for some $\lambda $ $%
> $ $2$ and some sequence of positive real numbers $\{\phi _{n}\}$, $\phi
_{n}$ $\geq $ $1$. Then $\sum_{i=1}^{p}\Theta _{i,n}^{(\cdot )}(\mathcal{E}%
_{n})$ $\simeq $ $\phi _{n}\mathcal{E}_{n}^{-\lambda }$ and $\Xi _{n}(%
\mathcal{E}_{n})$ $\simeq $ $\phi _{n}\mathcal{E}_{n}^{-(\lambda -1)}$.
Compare this to results in \citet[Corollary 2.2, Sect. 2.3]{ZhangCheng2018}
involving geometric decay. Similarly $\sum_{i=1}^{p}\check{\Theta}_{i,n}^{(%
\mathcal{M})(\cdot )}(\mathcal{E}_{n})$ $\simeq $ $\phi _{n}\mathcal{E}%
_{n}^{-\lambda }$ and $\check{\Xi}_{n}(\mathcal{E}_{n})$ $\simeq $ $\phi _{n}%
\mathcal{E}_{n}^{-(\lambda -1)}$. Then for any $\mathcal{E}_{n}$ $=$ $%
O(n^{3(1-1/q)/4})$ and $\gamma _{n}$ $\in $ $(0,1)$ it is straightforward to
verify that Assumption \ref{assume:ex_phys}.a,b hold. The second part of
Assumption \ref{assume:ex_phys}.c holds with $\xi _{n}$ $=$ $\phi _{n}$.

Now assume $q$ $=$ $2$, heterogeneity $\phi _{n}$ $=$ $n^{\phi }$, $\phi $ $<
$ $(3c$ $-$ $1/4)\wedge (1/2)$, memory decay $\lambda $ $=$ $2$, and set
block size $b_{n}$ $=$ $n^{b}$ for $b$ $\in $ $(0,\phi )$, and $\mathcal{E}%
_{n}$ $=$ $n^{c}$ and $\gamma _{n}$ $=$ $n^{-\gamma }$ for some $c$ $\in $ $%
(1/12,1/4)$ and $\gamma $ $>$ $0$. Thus $l_{n}$ $=$ $\ln (pn/\gamma _{n})$ $=
$ $\ln (pn^{1+\gamma })$ and $1\vee \ln (p/\mathcal{E}_{n}^{-\lambda })$ $%
\lesssim $ $\ln (pn)$ for $p$ $\geq $ $e$. If $p/n$ $\rightarrow $ $\infty $
then 
\begin{eqnarray*}
&&\rho _{n}\lesssim \left\{ n^{c/2-1/8}\ln (p)^{7/8}+n^{\phi /3-2c/3}\ln
\left( p\right) ^{2/3}\right\} +n^{-\gamma }+\left( \frac{1}{%
n^{3c/2-1/8-\phi /2}\ln (p)^{3/8}}\right) ^{2/3} \\
&&\rho _{n}^{\ast (\mathcal{M})}\lesssim n^{c/2-1/8}\ln (p)^{7/8}+n^{\phi
/3-2c/3}\ln \left( p\right) ^{2/3}+n^{-(\phi -b)/3}\ln (p)^{2/3} \\
&&\text{ \ \ \ \ \ \ \ \ \ \ \ \ \ \ \ }+n^{-\gamma }+\left( \frac{1}{%
n^{3c/2-1/8-\phi /2}\ln (p)^{3/8}}\right) ^{2/3}.
\end{eqnarray*}%
Hence from (\ref{Emax})%
\begin{equation*}
\mathbb{E}\max_{i}\bar{x}_{i,n}^{2}\lesssim \frac{2\ln \left( 2p\right) }{%
\mathcal{N}_{n}}\sqrt{\mathbb{E}\max_{i,l}\left( \frac{1}{b_{n}}\mathcal{S}%
_{n,l}^{(\mathcal{M})}(i)\right) ^{2}\mathbb{E}\left( \max_{i}\frac{1}{n}%
\sum_{l=1}^{\mathcal{N}_{n}}\left\vert \mathcal{S}_{n,l}^{(\mathcal{M}%
)}(i)\right\vert \right) ^{2}}+o(1)
\end{equation*}%
for any $p$ satisfying%
\begin{eqnarray}
\ln (p) &=&o\left( \frac{n^{3q/4}}{\max_{i,t}\left\Vert x_{i,t}\right\Vert
_{2r}^{r}\mathcal{M}_{n}^{3q/2-(r-1)}}\wedge \frac{n^{q/2-c/2+1/8}}{\mathcal{%
M}_{n}^{q}}\wedge \frac{n^{4q/7+4c/7-1/7}}{\mathcal{M}_{n}^{8q/7}}\right. 
\label{lnp} \\
&&\text{ \ \ \ \ \ \ \ }\left. \wedge \frac{\ln \left( \mathcal{M}%
_{n}^{q}/\max_{i}\mathbb{E}\left\vert \bar{x}_{i,n}\right\vert ^{q}\right) }{%
\left\Vert \max_{i,t}\left\vert x_{i,t}\right\vert \right\Vert
_{rq}^{qr/(r-1)}}\right) .  \notag
\end{eqnarray}%
Notice $3q/2$ $-$ $(r$ $-$ $1)$ $>$ $0$ for any $q$ $\geq $ $1$ when $r$ $<$ 
$2.5$. In that case the truncation point $\mathcal{M}_{n}$ exacts two forces
discussed in Remark \ref{rm:Mn}: larger $\mathcal{M}_{n}$ reduces the
truncation approximation error, allowing for larger $\ln (p)$; while smaller 
$\mathcal{M}_{n}$ reduces the tail remainder in the truncation
decomposition. Furthermore, \textit{ceteris paribus} if $c<q+1/4$, $\mathcal{%
M}_{n}$ $\simeq $ $n^{m}$, $m$ $>$ $0$ and%
\begin{equation*}
m<\frac{3q/4}{3q/2+r-1}\wedge \frac{q/2-c/2+1/8}{q}\wedge \frac{4q+4c-1}{8q},
\end{equation*}%
then the fourth term in (\ref{lnp}) dominates and larger $\mathcal{M}_{n}$
is optimal: $||\max_{i,t}\left\vert x_{i,t}\right\vert ||_{rq}^{qr/(r-1)}$ $%
\times $ $\ln (p)$ $=$ $o(m$ $\times $ $\ln \left( n/\max_{i}\mathbb{E}%
\left\vert \bar{x}_{i,n}\right\vert ^{q}\right) )$.

\section{Conclusion\label{sec:conc}}

We develop $\mathcal{L}_{q}$-maximal moment bounds similar to \cite%
{Nemirovski2000}, in a general dependence setting. Classic arguments exploit
symmetrization under independence. Under arbitrary weak dependence we
exploit a standard multiplier blocking argument with a negligible truncation
approximation, and use Gaussian approximation and comparison theory to
sidestep symmetrization. We also require a concentration bound for tail
probabilities that appears new, and works roughly like a Nemirovski bound
for tail measures. The latter arises from the truncation approximation
error, leading to a higher moment requirement than in Nemirovski's case with
independence. Examples are provided in order to verify assumptions and to
yield Gaussian approximations, working under either mixing or physical
dependence conditions. We do not focus on sharpness, leaving that idea for
future research, and instead seek broad conditions such that an $\mathcal{L}%
_{q}$-maximal moment is bounded.

\section{Appendix: technical proofs}

\label{app:proofs}

Recall%
\begin{eqnarray*}
&&\rho _{n}:=\sup_{z\geq 0}\left\vert \mathbb{P}\left( \max_{i}\left\vert 
\mathcal{X}_{n}(i)\right\vert \leq z\right) -\mathbb{P}\left(
\max_{i}\left\vert \boldsymbol{X}_{n}(i)\right\vert \leq z\right)
\right\vert  \\
&&\rho _{n}^{\ast (\mathcal{M})}:=\sup_{z\geq 0}\left\vert \mathbb{P}\left(
\max_{i}\left\vert \mathcal{X}_{n}^{(\mathcal{M})\ast }(i)\right\vert \leq
z\right) -\mathbb{P}\left( \max_{i}\left\vert \boldsymbol{X}_{n}^{(\mathcal{M%
})}(i)\right\vert \leq z\right) \right\vert  \\
&&\delta _{n}^{(\mathcal{M})}:=\sup_{z\geq 0}\left\vert \mathbb{P}\left(
\max_{i}\left\vert \boldsymbol{X}_{n}^{(\mathcal{M})}(i)\right\vert \leq
z\right) -\mathbb{P}\left( \max_{i}\left\vert \boldsymbol{X}%
_{n}(i)\right\vert \leq z\right) \right\vert .
\end{eqnarray*}%
\noindent \textbf{Proof of Proposition \ref{prop:main_unbound}.} By a change
of variables,%
\begin{eqnarray*}
\mathbb{E}\max_{i}\left\vert \bar{x}_{i,n}\right\vert ^{q} &=&\mathbb{E}%
\left( \max_{i}\left\vert \bar{x}_{i,n}\right\vert ^{q}\mathcal{I}_{\max_{i}|%
\bar{x}_{i,n}|\leq \mathcal{M}_{n}}\right) +\mathbb{E}\left(
\max_{i}\left\vert \bar{x}_{i,n}\right\vert ^{q}\mathcal{I}_{\max_{i}|\bar{x}%
_{i,n}|>\mathcal{M}_{n}}\right)  \\
&=&q\int_{0}^{\mathcal{M}_{n}}u^{q-1}\mathbb{P}\left( \max_{i}\left\vert 
\bar{x}_{i,n}\right\vert >u\right) du+\mathbb{E}\left( \max_{i}\left\vert 
\bar{x}_{i,n}\right\vert ^{q}\mathcal{I}_{\max_{i}|\bar{x}_{i,n}|>\mathcal{M}%
_{n}}\right)  \\
&=&\mathfrak{I}_{n,1}+\mathfrak{I}_{n,2}.
\end{eqnarray*}%
\textbf{Step 1 (}$\mathfrak{I}_{n,1}$\textbf{)}. Add and subtract like
terms, and invoke the triangle inequality and the proof of Proposition \ref%
{prop:main} to yield for some $r$ $>$ $1$ such that $%
\max_{i,t}||x_{i,t}||_{2r}$ $<$ $\infty $, and any $x$ $\geq $ $0$, 
\begin{eqnarray}
&&\left\vert \mathbb{P}\left( \sqrt{n}\max_{i}\left\vert \bar{x}%
_{i,n}\right\vert \leq x\right) -\mathbb{P}\left( \max_{i}\left\vert \frac{1%
}{\sqrt{n}}\sum_{l=1}^{\mathcal{N}_{n}}\varepsilon _{l}\mathcal{S}_{n,l}^{(%
\mathcal{M})}(i)\right\vert \leq x\right) \right\vert   \label{PPM} \\
&&\text{ \ \ \ \ \ \ \ \ \ \ \ }\leq \delta _{n}^{(\mathcal{M})}+\left\{
\rho _{n}+\rho _{n}^{\ast (\mathcal{M})}\right\}   \notag \\
&&\text{ \ \ \ \ \ \ \ \ \ \ \ }\leq \left\{ \left( \frac{n}{\mathcal{M}%
_{n}^{r-1}}\max_{i,t}\left\Vert x_{i,t}\right\Vert _{2r}^{2r}\right)
^{1/2}\ln (p)\right\} ^{2/3}+\left\{ \rho _{n}+\rho _{n}^{\ast (\mathcal{M}%
)}\right\} :=\mathcal{B}_{n}(p)+\mathcal{C}_{n}(p).  \notag
\end{eqnarray}%
Therefore, by twice a change of variables%
\begin{eqnarray*}
\mathfrak{I}_{n,1} &=&q\int_{0}^{\mathcal{M}_{n}}u^{q-1}\mathbb{P}\left( 
\sqrt{n}\max_{i}\left\vert \bar{x}_{i,n}\right\vert >\sqrt{n}u\right) du \\
&=&\frac{q}{n^{q/2}}\int_{0}^{\mathcal{M}_{n}}v^{q-1}\mathbb{P}\left( \sqrt{n%
}\max_{i}\left\vert \bar{x}_{i,n}\right\vert >v\right) dv \\
&\leq &\frac{q}{n^{q/2}}\int_{0}^{\mathcal{M}_{n}}v^{q-1}\mathbb{P}\left(
\max_{i}\left\vert \frac{1}{\sqrt{n}}\sum_{l=1}^{\mathcal{N}_{n}}\varepsilon
_{l}\mathcal{S}_{n,l}^{(\mathcal{M})}(i)\right\vert \leq v\right) dv+\frac{%
q\left\{ \mathcal{B}_{n}+\mathcal{C}_{n}(p)\right\} }{n^{q/2}}\int_{0}^{%
\mathcal{M}_{n}}v^{q-1}dv \\
&=&\mathbb{E}\max_{i}\left\vert \frac{1}{n}\sum_{l=1}^{\mathcal{N}%
_{n}}\varepsilon _{l}\mathcal{S}_{n,l}^{(\mathcal{M})}(i)\right\vert ^{q}+%
\frac{\mathcal{M}_{n}^{q}}{n^{q/2}}\left\{ \mathcal{B}_{n}(p)+\mathcal{C}%
_{n}(p)\right\} .
\end{eqnarray*}%
Replicate (\ref{EEx})-(\ref{Esum*1}), and use $n$ $=$ $\mathcal{N}_{n}b_{n}$%
\ to complete the bounds:%
\begin{eqnarray*}
\mathfrak{I}_{n,1} &\leq &\left( \frac{2\ln \left( 2p\right) }{\mathcal{N}%
_{n}}\right) ^{q/2}\sqrt{\mathbb{E}\max_{i,l}\left\vert \frac{1}{b_{n}}%
\mathcal{S}_{n,l}^{(\mathcal{M})}(i)\right\vert ^{q}\mathbb{E}\left( \max_{i}%
\frac{1}{\mathcal{N}_{n}}\sum_{l=1}^{\mathcal{N}_{n}}\left\vert \mathcal{S}%
_{n,l}^{(\mathcal{M})}(i)\right\vert \right) ^{q}} \\
&&\text{ \ \ \ \ \ \ \ }+\frac{\mathcal{M}_{n}^{q}}{n^{q/2}}\left\{ \mathcal{%
B}_{n}(p)+\mathcal{C}_{n}(p)\right\} .
\end{eqnarray*}%
\textbf{Step 2 (}$\mathfrak{I}_{n,2}$\textbf{)}: Recall $\mathbb{\bar{P}}_{%
\mathcal{M}_{n}}$ $:=$ $\max_{i}\mathbb{P}(|\bar{x}_{i,n}|$ $>$ $\mathcal{M}%
_{n})$. H\"{o}lder's inequality yields for $r$ $>$ $1$%
\begin{equation}
\mathfrak{I}_{n,2}=\mathbb{E}\max_{i}\left\vert \bar{x}_{i,n}\right\vert ^{q}%
\mathcal{I}_{\max_{i}|\bar{x}_{i,n}|>\mathcal{M}_{n}}\leq \left\Vert
\max_{i}\left\vert \bar{x}_{i,n}\right\vert \right\Vert _{rq}^{q}\times 
\mathbb{\bar{P}}_{\mathcal{M}_{n}}^{(r-1)/r}.  \label{Jn2}
\end{equation}%
\textbf{Claim (a). }Combine Lemma \ref{lm:concen}.a with (\ref{Jn2}) and
Minkowski's inequality to deduce $\mathfrak{I}_{n,2}$ $\leq $ $%
2^{(r-1)/r}||\max_{i}\left\vert \bar{x}_{i,n}\right\vert ||_{rq}^{q}$ $%
\times $ $\{\ln (p)/\ln \mathbb{\bar{P}}_{\mathcal{M}_{n}}^{-1}\}^{(r-1)/r}$%
. Now use Lyapunov's inequality and $\max_{i}\left\vert a_{i}\right\vert $ $%
\leq $ $\sum \left\vert a_{i}\right\vert $ to obtain for $s$ $>$ $q$,%
\begin{equation*}
\mathfrak{I}_{n,2}\leq 2^{(r-1)/r}p^{q/(sr)}\max_{i}\left\Vert \bar{x}%
_{i,n}\right\Vert _{rs}^{q}\left( \frac{\ln (p)}{\ln \mathbb{\bar{P}}_{%
\mathcal{M}_{n}}^{-1}}\right) ^{(r-1)/r}:=\mathcal{D}_{n}(p).
\end{equation*}%
\textbf{Claim (b).} By Lemma \ref{lm:concen}.b with (\ref{Jn2}), and
Lyapunov's inequality, for $s$ $>$ $q$ and $r$ $>$ $1$%
\begin{equation*}
\mathfrak{I}_{n,2}\leq 2^{(r-1)/r}p^{q/(sr)}p^{q/(sr)}\max_{i}\left\Vert 
\bar{x}_{i,n}\right\Vert _{rs}^{q}\left( \frac{\ln (p)}{\ln \left( \mathcal{M%
}_{n}^{q}/\max_{i}\mathbb{E}\left\vert \bar{x}_{i,n}\right\vert ^{q}\right) }%
\right) ^{(r-1)/r}:=\mathcal{D}_{n}(p).
\end{equation*}%
Moreover, if $\mathbb{E}|\bar{x}_{i,n}|^{s}$ $=$ $O(n^{-a})$ for some $a$ $>$
$0$ and $s$ $>$ $1$ then by Lyapunov's inequality and $\max_{i}|z_{i}|$ $%
\leq $ $\sum_{i=1}^{p}|z_{i}|$ we have for $s$ $>$ $rq$%
\begin{equation*}
\left\Vert \max_{i}\left\vert \bar{x}_{i,n}\right\vert \right\Vert
_{rq}^{q}\leq p^{q/s}\max_{i}\left\Vert \bar{x}_{i,n}\right\Vert
_{s}^{q}\leq K\frac{p^{q/s}}{n^{aq/s}}.
\end{equation*}%
Hence, in view of (\ref{Jn2}) we may replace $||\max_{i,t}\left\vert
x_{i,t}\right\vert ||_{rq}^{q}$ with $Kp^{q/s}/n^{aq/s}$ yielding as claimed 
$\mathfrak{I}_{n,2}$ $\leq $ $2^{(r-1)/r}p^{q/s}n^{-aq/s}\left( \ln (p)/\ln (%
\mathcal{M}_{n}^{q}n^{a})\right) ^{(r-1)/r}$.\medskip \newline
\textbf{Claim (c). }Let $\max_{i}\mathbb{P}(|\bar{x}_{i,n}|$ $\geq c)$ $\leq 
$ $a\exp \{-bn^{\gamma }c^{\gamma }\}$, $a,b$ $>$ $0$, $\gamma $ $>$ $q$. By
Jensen's inequality 
\begin{equation*}
\mathfrak{I}_{n,2}\leq \frac{1}{\lambda }\ln \left( p\max_{i}\mathbb{E}%
\left( \exp \left\{ \lambda \left\vert \bar{x}_{i,n}\right\vert ^{q}\mathcal{%
I}_{\max_{i}|\bar{x}_{i,n}|>\mathcal{M}_{n}}\right\} \right) \right) \text{
for any }\lambda >0.
\end{equation*}%
Twice a change of variables, sub-exponential tails, and $\mathcal{M}_{n}$ $%
\rightarrow $ $\infty $ yield for any $c$ $\in $ $[q.\gamma )$, and any $b$ $%
\in $ $(0,a]$ 
\begin{eqnarray*}
&&\mathbb{E}\left( \exp \left\{ \lambda \left\vert \bar{x}_{i,n}\right\vert
^{q}\mathcal{I}_{\max_{i}|\bar{x}_{i,n}|>\mathcal{M}_{n}}\right\} \right)  \\
&&\text{ \ \ \ \ \ \ \ }=1+\int_{\exp \left\{ \lambda \mathcal{M}%
_{n}^{q}\right\} }^{\infty }\mathbb{P}\left( \left\vert \bar{x}%
_{i,n}\right\vert >\frac{1}{\lambda ^{1/q}}\ln (u)^{1/q}\right) du \\
&&\text{ \ \ \ \ \ \ \ }=1+q\lambda \int_{\exp \left\{ \lambda \mathcal{M}%
_{n}^{q}\right\} }^{\infty }\mathbb{P}\left( \left\vert \bar{x}%
_{i,n}\right\vert >v\right) \exp \{\lambda v^{q}\}v^{q-1}dv \\
&&\text{ \ \ \ \ \ \ \ }\leq 1+2q\lambda \int_{\exp \left\{ \lambda \mathcal{%
M}_{n}^{q}\right\} }^{\infty }\exp \{\lambda v^{q}-bn^{\gamma }v^{\gamma
}\}v^{q-1}dv \\
&&\text{ \ \ \ \ \ \ \ }\leq 1+2q\lambda \int_{\exp \left\{ \lambda \mathcal{%
M}_{n}^{q}\right\} }^{\infty }\exp \left\{ -bn^{\gamma }v^{c}\right\} dv%
\text{ \ }\forall n\geq \bar{n}\text{ and some }\bar{n}\in \mathbb{N} \\
&&\text{ \ \ \ \ \ \ \ }=1+\frac{2q\lambda }{cn^{\gamma /c}b^{1/c}}%
\int_{\exp \left\{ \lambda \mathcal{M}_{n}^{q}\right\} }^{\infty }\frac{1}{%
u^{1/c-1}}\exp \left\{ -u\right\} du \\
&&\text{ \ \ \ \ \ \ \ }\leq 1+\frac{2q\lambda }{cn^{\gamma /c}b^{1/c}}%
\int_{\exp \left\{ \lambda \mathcal{M}_{n}^{q}\right\} }^{\infty }\exp
\left\{ -u\right\} du=1+\frac{2q\lambda }{cn^{\gamma /c}b^{1/c}}\exp \left\{
-\exp \left\{ \lambda \mathcal{M}_{n}^{q}\right\} \right\} .
\end{eqnarray*}%
Now exploit $\ln (1$ $+$ $z)$ $\leq $ $z$ for $z$ $>$ $0$ and $\lambda $ $%
\geq $ $1$ to yield%
\begin{eqnarray*}
\mathfrak{I}_{n} &\leq &\frac{1}{\lambda }\ln \left( p\left[ 1+\frac{%
2q\lambda }{cn^{\gamma /c}b^{1/c}}\exp \left\{ -\exp \left\{ \lambda 
\mathcal{M}_{n}^{q}\right\} \right\} \right] \right)  \\
&\leq &\frac{1}{\lambda }\ln (p)+\frac{2q\lambda }{cn^{\gamma /c}b^{1/c}}%
\exp \left\{ -\exp \left\{ \lambda \mathcal{M}_{n}^{q}\right\} \right\} \leq 
\frac{1}{\lambda }\ln (p)+\frac{2q\lambda }{cn^{\gamma /c}b^{1/c}}\exp
\left\{ -\exp \left\{ \mathcal{M}_{n}^{q}\right\} \right\} \text{.}
\end{eqnarray*}%
Choose $\lambda $ $=$ $\sqrt{cn^{\gamma /c}b^{1/c}/(2q)}\times \exp \{.5\exp
\left\{ \mathcal{M}_{n}^{q}\right\} \}\sqrt{\ln (p)}$\ to minimize the upper
bound. Thus $\lambda $ $\geq $ $1$ when $p$ $\geq $ $e$ and $\mathcal{M}_{n}$
$\geq $ $[\ln (1$ $\vee $ $\ln (2q/\{cb^{1/c}\}))]^{1/q}$ , while the latter
holds for all $n\geq $ \b{n} and some \b{n} $\in $ $\mathbb{N}$. Thus%
\begin{equation*}
\mathfrak{I}_{n,2}\leq 2\sqrt{\frac{2q}{cn^{\gamma /c}b^{1/c}}}\frac{\sqrt{%
\ln (p)}}{\exp \left\{ \exp \left\{ \mathcal{M}_{n}^{q}\right\} /2\right\} }%
:=\mathcal{D}_{n}(p).
\end{equation*}%
Now choose $c$ $=$ $q$ and $b$ $=$ $a$ to complete the proof.\footnote{%
If, for example, $\mathcal{M}_{n}$ $=$ $n^{m}$, $m$ $>0$, then $\lambda $ $%
\geq $ $1$ for all $n$ $\geq $ $1$ $\vee $ $[\ln (1\vee \ln
(2/a^{1/q}))]^{1/(qm)}$ $:=$ \b{n}.
\par
{}} $\mathcal{QED}$.\medskip \newline
\textbf{Proof of Lemma \ref{lm:concen}. }Define $\mathbb{P}_{\mathcal{M}_{n}}
$ $:=$ $\max_{i}\mathbb{P}(|\bar{x}_{i,n}|$ $\leq $ $\mathcal{M}_{n})$ and $%
\mathbb{\bar{P}}_{\mathcal{M}_{n}}$ $:=$ $\max_{i}\mathbb{P}(|\bar{x}_{i,n}|$
$>$ $\mathcal{M}_{n})$. We first prove for any $\lambda $ $>$ $0$ 
\begin{equation}
\mathbb{P}\left( \max_{i}\left\vert \bar{x}_{i,n}\right\vert \geq \mathcal{M}%
_{n}\right) \leq \frac{1}{\lambda }\ln (p)+\frac{1}{\lambda }\exp \left\{
\lambda \right\} \times \mathbb{\bar{P}}_{\mathcal{M}_{n}}.  \label{PP}
\end{equation}%
Jensen's inequality yields for any $\lambda $ $>$ $0$ 
\begin{eqnarray}
\mathbb{P}\left( \max_{i}\left\vert \bar{x}_{i,n}\right\vert >\mathcal{M}%
_{n}\right)  &=&\frac{1}{\lambda }\mathbb{E}\left[ \ln \left( \exp \left\{
\lambda \mathcal{I}_{\max_{i}|\bar{x}_{i,n}|\geq \mathcal{M}_{n}}\right\}
\right) \right]   \notag \\
&\leq &\frac{1}{\lambda }\ln \left( \mathbb{E}\left[ \exp \left\{ \lambda 
\mathcal{I}_{\max_{i}|\bar{x}_{i,n}|\geq \mathcal{M}_{n}}\right\} \right]
\right)   \notag \\
&\leq &\frac{1}{\lambda }\ln \left( p\max_{i}\mathbb{E}\left[ \exp \left\{
\lambda \mathcal{I}_{|\bar{x}_{i,n}|\geq \mathcal{M}_{n}}\right\} \right]
\right) .  \notag
\end{eqnarray}%
By construction $\max_{i}\mathbb{E}[\exp \{\lambda \mathcal{I}_{|\bar{x}%
_{i,n}|\geq \mathcal{M}_{n}}\}]$ $\leq $ $\exp \{\lambda \}\mathbb{\bar{P}}_{%
\mathcal{M}_{n}}$ $+$ $\mathbb{P}_{\mathcal{M}_{n}}$ $\leq $ $\exp \{\lambda
\}\mathbb{\bar{P}}_{\mathcal{M}_{n}}$ $+$ $1$. Now $\ln (1$ $+$ $x)$ $\leq $ 
$x$ $\forall x$ $\geq $ $0$ achieves (\ref{PP}).\medskip \newline
\textbf{Claim (a). }The upper bound in (\ref{PP}) is minimized at $\lambda $ 
$=$ $\ln \mathbb{\bar{P}}_{\mathcal{M}_{n}}^{-1}$ $+$ $\ln \ln (p)$ as $n$ $%
\rightarrow $ $\infty $ Thus%
\begin{eqnarray*}
\mathbb{P}\left( \max_{i}\left\vert \bar{x}_{i,n}\right\vert \geq \mathcal{M}%
_{n}\right)  &\leq &\frac{1}{\ln \left( \mathbb{\bar{P}}_{\mathcal{M}%
_{n}}^{-1}\ln (p)\right) }\ln (p)+\frac{1}{\ln \left( \mathbb{\bar{P}}_{%
\mathcal{M}_{n}}^{-1}\ln (p)\right) }\exp \left\{ \ln \left( \mathbb{\bar{P}}%
_{\mathcal{M}_{n}}^{-1}\ln (p)\right) \right\} \times \mathbb{\bar{P}}_{%
\mathcal{M}_{n}} \\
&=&2\frac{\ln (p)}{\ln \mathbb{\bar{P}}_{\mathcal{M}_{n}}^{-1}+\ln \ln (p)}%
\leq 2\frac{\ln (p)}{\ln \mathbb{\bar{P}}_{\mathcal{M}_{n}}^{-1}},
\end{eqnarray*}%
proving the claim when $p$ $>$ $e$ and $p$ $=$ $o(\mathbb{\bar{P}}_{\mathcal{%
M}_{n}}^{-1})$, given $\mathbb{\bar{P}}_{\mathcal{M}_{n}}^{-1}$ $>$ $1$%
.\medskip \newline
\textbf{Claim (b).} Use $\mathbb{\bar{P}}_{\mathcal{M}_{n}}$ $\leq $ $%
\mathcal{M}_{n}^{-q}\max_{i}\mathbb{E}|\bar{x}_{i,n}|^{q}$ with (\ref{PP})
and $\mathcal{L}_{q}$-boundedness. Minimizing the upper bound with respect
to $\lambda $, and setting $p$ $=$ $o(\mathcal{M}_{n}^{q}/\max_{i}\mathbb{E}%
\left\vert \bar{x}_{i,n}\right\vert ^{q})$ yields the claim.\medskip \newline
\textbf{Claim (c). }Let $\mathbb{\bar{P}}_{\mathcal{M}_{n}}$ $\leq $ $a\exp
\{-bn^{\gamma }\mathcal{M}_{n}^{\gamma }\}$. By (\ref{PP}) with $\lambda $ $=
$ $n^{\phi }\mathcal{M}_{n}^{\phi }\ln (\ln p)$ for any $\phi $ $\in $ $%
(0,\gamma )$ we have 
\begin{eqnarray*}
\mathbb{P}\left( \max_{i}\left\vert \bar{x}_{i,n}\right\vert \geq \mathcal{M}%
_{n}\right)  &\leq &\frac{1}{\lambda }\ln (p)+\frac{1}{\lambda }\exp \left\{
\lambda \right\} a\exp \left\{ -bn^{\gamma }\mathcal{M}_{n}^{\gamma
}\right\}  \\
&=&\frac{\ln (p)}{n^{\phi }\mathcal{M}_{n}^{\phi }\ln (\ln p)}+a\frac{\left(
\ln (p)\right) ^{n^{\phi }\mathcal{M}_{n}^{\phi }}}{n^{\phi }\mathcal{M}%
_{n}^{\phi }\exp \left\{ bn^{\gamma }\mathcal{M}_{n}^{\gamma }\right\} \ln
(\ln p)}.
\end{eqnarray*}%
The upper bound is $o(1)$ if $p$ $>$ $e$ and $\ln (p)$ $=$ $o(n^{\phi }%
\mathcal{M}_{n}^{\phi }\wedge \exp \{bn^{\gamma -\phi }\mathcal{M}%
_{n}^{\gamma -\phi }\})$ $=$ $o(n^{\phi }\mathcal{M}_{n}^{\phi })$. $%
\mathcal{QED}$.\medskip \newline
\textbf{Proof of Lemma \ref{lm:ex1}.} Conditions 1-3 in \cite%
{ChangChenWu2024} hold by construction. In particular, their Condition 2 $%
\alpha _{i,j,n}(m)$ $\leq $ $a\exp \{-bl^{\gamma }\}$ holds for $a$ $\geq $ $%
1$, some $b$ $>$ $0$ and $\gamma $ $=$ $1$. Thus $\rho _{n}$ $\lesssim $ $%
\chi _{n}^{2/3}\ln (p)/n^{1/9}$ $+$ $\chi _{n}\left( \ln (p)\right)
^{7/6}/n^{1/9}$ by their Theorem 1 and the mapping theorem, yielding (\ref%
{rex1}).

Now turn to $\mathcal{X}_{n}^{\ast }(i)$ $=$ $1/\sqrt{n}\sum_{t=1}^{n}\eta
_{t}x_{i,t}$, and let $\{\boldsymbol{X}_{n}^{\ast }(i)\}$ be a Gaussian
process with marginals $\boldsymbol{X}_{n}^{\ast }(i)$ $\sim $ $N(0,\mathbb{E%
}\mathcal{X}_{n}^{\ast 2}(i))$. We will derive the following Gaussian
approximation and Gaussian (without blocking)-to-Gaussian (with blocking)
comparison bounds for some $\{c_{n}(p),d_{n}(p)\}$, 
\begin{eqnarray}
&&\rho _{n}^{\ast \ast }:=\sup_{z\geq 0}\left\vert \mathbb{P}\left(
\max_{i}\left\vert \mathcal{X}_{n}^{\ast }(i)\right\vert \leq z\right) -%
\mathbb{P}\left( \max_{i}\left\vert \boldsymbol{X}_{n}^{\ast }(i)\right\vert
\leq z\right) \right\vert \lesssim c_{n}(p)  \label{rho*_til} \\
&&\delta _{n}^{\ast }:=\sup_{z\geq 0}\left\vert \mathbb{P}\left(
\max_{i}\left\vert \boldsymbol{X}_{n}(i)\right\vert \leq z\right) -\mathbb{P}%
\left( \max_{i}\left\vert \boldsymbol{X}_{n}^{\ast }(i)\right\vert \leq
z\right) \right\vert \lesssim d_{n}(p).  \label{d*}
\end{eqnarray}%
Hence by the triangle inequality%
\begin{equation*}
\rho _{n}^{\ast }=\sup_{z\geq 0}\left\vert \mathbb{P}\left(
\max_{i}\left\vert \mathcal{X}_{n}^{\ast }(i)\right\vert \leq z\right) -%
\mathbb{P}\left( \max_{i}\left\vert \boldsymbol{X}_{n}(i)\right\vert \leq
z\right) \right\vert \lesssim c_{n}(p)\vee d_{n}(p).
\end{equation*}%
\textbf{Step 1: eq. (\ref{rho*_til}).} $\{\eta _{t}x_{i,t}\}_{t=1}^{n}$ is $%
\alpha $-mixing with coefficients $\alpha _{i,j,n}((m$ $-$ $b_{n})_{+})$ by $%
b_{n}$-dependence of $\eta _{t}$, mutual independence and measurability. If
the sub-exponential tail component $\chi _{n}$\ is a fixed constant$\chi _{n}
$ $=$ $K$ then we can use Proposition 3 in \cite{ChangJiangShao2023}.

Otherwise, we use the fact that uniform geometric $\alpha $-mixing implies
geometric physical dependence in the sense of \cite{Wu2005} and \cite%
{WuMin2005}. Recall $x_{i,t}$ $=$ $g_{i,t}(\epsilon _{i,t},$ $\epsilon
_{i,t-1},\ldots )$, where $\epsilon _{i,t}$ are iid for each $i$. Let $%
\{\epsilon _{i,t}^{\prime }\}$ be an independent copy of $\{\epsilon _{i,t}\}
$, and define the coupled process $x_{i,t}^{\prime }(m)$ $:=$ $%
g_{i,t}(\epsilon _{i,t},\ldots ,\epsilon _{i,t-m+1},$ $\epsilon
_{i,t-m}^{\prime },$ $\epsilon _{i,t-m-1},\ldots )$, $m$ $=$ $0,1,2,...$ Now
define the $\mathcal{L}_{p}$-\textit{physical dependence} measure $\theta
_{i,t}^{(p)}(m)$ $:=$ $||x_{i,t}-x_{i,t}^{\prime }(m)||_{p}$. We have $%
\theta _{i,n,t}^{(q)}(m)$ $\leq $ $2^{1+2/q}||x_{i,t}||_{q}^{1/q}\alpha
_{i,j,n}^{1/(qr)}(m)$ for $1/q$ $+$ $1/r$ $=$ $1$: see \citet[Lemma
1.6]{McLeish1975} and \citet[Theorem 2.1]{Hill2025_mixg}, cf. \citet[proof
of Theorem 2.10]{Hill2024_maxlln}. Hence by supposition and therefore moment
bound (\ref{Exq}), $q$ $\geq $ $2$, and $r$ $=$ $q/(q$ $-$ $1)$, we have for
some $\omega $ $\in $ $(0,1)$%
\begin{equation*}
\max_{i,t}\theta _{i,n,t}^{(q)}(m)\leq 2^{1+2/q}q^{1/q}\chi _{n}^{1/q}\omega
^{m/(qr)}\leq 4q^{1/q}\chi _{n}^{1/q}\omega ^{m(q-1)/q^{2}}.
\end{equation*}

Moreover, by mutual independence and additivity $y_{i,t}$ $:=$ $\eta
_{t}x_{i,t}$ is also geometrically $\mathcal{L}_{q/2}$-physical dependent.
Let $\{\varepsilon _{l}^{(j)}\}_{l=1}^{\mathcal{N}_{n}}$ be independent
copies of $\{\varepsilon _{l}\}_{l=1}^{\mathcal{N}_{n}}$, $j$ $=$ $%
2,...,b_{n}$, and let $\varepsilon _{l}^{(1)}$ $=$ $\varepsilon _{l}$. Let $%
\{\tilde{\varepsilon}_{t}\}_{t=1}^{n}$ be iid random variables constructed
block-wise $\tilde{\varepsilon}_{b_{n}(l-1)+j}$ $=$ $\varepsilon _{l}^{(j)}$
for $j$ $=$ $1,...,b_{n}$. Hence $y_{i,t}=h(\tilde{\epsilon}_{i,t},\tilde{%
\epsilon}_{i,t-1},\ldots )$ where $\tilde{\epsilon}_{i,t}$ $=$ $[\epsilon
_{i,t},\tilde{\varepsilon}_{t}]$. Using the above notation define $\check{%
\theta}_{i,n,t}^{(q)}(m)$ $:=$ $||\eta _{t}$ $-$ $\eta _{t}^{\prime
}(m)||_{q}$ $\leq $ $\mathcal{I}_{m\leq b_{n}}$ and $\tilde{\theta}%
_{i,n,t}^{(r)}(m)$ $:=$ $||y_{i,t}$ $-$ $y_{i,t}^{\prime }(m)||_{r}$. Then
by Minkowski and Cauchy-Schwartz inequalities, $|\eta _{t}|$ $\leq $ $c$ $%
a.s.$ for some $c$ $\in $ $(0,\infty )$, $\chi _{n}$ $\geq $ $1$\ and (\ref%
{Exq}), 
\begin{eqnarray}
\tilde{\theta}_{i,n,t}^{(q/2)}(m) &\leq &\left\Vert x_{i,t}\right\Vert
_{q}\left\Vert \eta _{t}-\eta _{t}^{\prime }(m)\right\Vert _{q}+\left\Vert
\eta _{t}^{\prime }(m)\right\Vert _{q}\left\Vert x_{i,t}-x_{i,t}^{\prime
}(m)\right\Vert _{q}  \notag \\
&=&\left\Vert x_{i,t}\right\Vert _{q}\check{\theta}_{i,n,t}^{(q)}(m)+\left%
\Vert \eta _{t}^{\prime }(m)\right\Vert _{q}\theta _{i,n,t}^{(q)}(m)  \notag
\\
&\lesssim &q\chi _{n}\mathcal{I}_{m\leq b_{n}}+4cq^{1/q}\chi
_{n}^{1/q}\omega ^{m(q-1)/q^{2}}\leq 4q\chi _{n}\left( \mathcal{I}_{m\leq
b_{n}}+c\omega ^{m(q-1)/q^{2}}\right) .  \label{theta_til1}
\end{eqnarray}%
Moreover,\ by mutual independence and the definition of $\chi _{n}$, $||\eta
_{t}x_{i,t}||_{\psi _{\zeta }}$ $\leq $ $c\chi _{n}$. Finally, by
supposition and the mixing assumption $\mathbb{E}\mathcal{X}_{n}^{\ast 2}(i)$
$\in $ $[K,\infty )$, cf. \cite{Davydov1968}. Theorem 3(\textit{ii}) in \cite%
{ChangChenWu2024} thus applies: for some $\beta ,\nu $ $\in $ $(0,\infty )$ 
\begin{equation*}
\rho _{n}^{\ast }\lesssim \frac{\chi _{n}\left( \ln (p)\right) ^{7/6}+\left(
\Psi _{n,\beta }^{(2)}\right) ^{1/3}\left( \Psi _{n,0}^{(2)}\right)
^{1/3}\left( \ln (p)\right) ^{2/3}}{n^{\beta /(12+6\beta )}}+\frac{\Phi
_{n,\beta ,\nu }\left( \ln (p)\right) ^{1+\nu }}{n^{\beta /(4+\beta )}},
\end{equation*}%
where $(\Psi _{n,\beta }^{(q)},\Phi _{n,\beta ,\nu })$ are \textit{%
aggregated dependence adjusted norms} (cf. \cite{WuWu2016,ChangJiangShao2023}%
):%
\begin{eqnarray*}
\Psi _{n,\beta }^{(q)} &=&\max_{i}\left\{ \sup_{m\geq 0}\left\{ \left(
m+1\right) ^{\beta }\sum_{i=m}^{\infty }\tilde{\theta}_{i,n,t}^{(q)}(m)%
\right\} \right\} \\
&\leq &4q\chi _{n}\left( \sup_{m\geq 0}\left\{ \left( m+1\right) ^{\beta
}\sum_{i=m}^{\infty }\mathcal{I}_{m\leq b_{n}}\right\} +\sup_{m\geq
0}\left\{ \left( m+1\right) ^{\beta }\sum_{i=m}^{\theta }\omega
^{l(q-1)/q^{2}}\right\} \right) \\
&\lesssim &4q\chi _{n}\left\{ b_{n}^{\beta }+\sup_{m\geq 0}\frac{\left(
m+1\right) ^{\beta }\omega ^{m(q-1)/q^{2}}}{1-\omega ^{(q-1)/q^{2}}}\right\}
\\
&\lesssim &4q\chi _{n}\left\{ b_{n}^{\beta }+\frac{\frac{q^{2\beta }}{\left(
q-1\right) ^{\beta }}\omega ^{\beta /\ln \left( 1/\omega \right)
-(q-1)/q^{2}}}{1-\omega ^{(q-1)/q^{2}}}\right\} \lesssim \chi
_{n}b_{n}^{\beta }q^{\beta +1}
\end{eqnarray*}%
and%
\begin{eqnarray}
\Phi _{n,\beta ,\nu } &=&\sup_{q\geq 2}\left\{ q^{-\nu }\sup_{m\geq
0}\left\{ \left( m+1\right) ^{\beta }\sum_{l=m}^{\infty }\tilde{\theta}%
_{i,n,t}^{(q)}(l)\right\} \right\}  \label{Phi} \\
&\lesssim &\sup_{q\geq 2}\left\{ q^{-\nu }\chi _{n}b_{n}^{\beta }q^{\beta
+1}\right\} =\chi _{n}b_{n}^{\beta }\text{ for }\nu =\beta +1,  \notag
\end{eqnarray}%
where $\omega $ $\in $ $(0,1)$, $q$ $\geq $ $2$, and $b_{n}$ $\rightarrow $ $%
\infty $ yield the upper bounds. Thus using $\nu $ $=$ $\beta $ $+$ $1$ we
have shown for some $\beta $ $>$ $0$ (which may be arbitrarily large under
geometric mixing) 
\begin{eqnarray}
\rho _{n}^{\ast } &\lesssim &\frac{\chi _{n}\left( \ln (p)\right)
^{7/6}+\chi _{n}^{2/3}b_{n}^{2\beta /3}\left( \ln (p)\right) ^{2/3}}{%
n^{\beta /(12+6\beta )}}+\frac{\chi _{n}b_{n}^{\beta }\left( \ln (p)\right)
^{\beta +2}}{n^{\beta /(4+\beta )}}  \label{rhotil} \\
&\leq &\frac{\chi _{n}b_{n}^{\beta }\left( \ln (p)\right) ^{\beta +2}}{%
n^{\beta /(4+\beta )}}:=c_{n}(p).  \notag
\end{eqnarray}%
Finally, the upper bound is minimized with $\beta ^{\ast }$ $=$ $\sqrt{4\ln
n/\{\ln b_{n}+\ln \ln (p)\}}$ $-$ $4$, where $\beta ^{\ast }$ $>$ $0$ given $%
\ln (p)$ $<$ $n^{1/4}/b_{n}$.\medskip \newline
\textbf{Step 2: eq. (\ref{d*}).} Define $\sigma _{n}^{2}(i,j)$ $:=$ $\mathbb{%
E}\mathcal{X}_{n}(i)\mathcal{X}_{n}(j)$,\ $\sigma _{n}^{\ast 2}(i,j)$ $:=$ $%
\mathbb{E}\mathcal{X}_{n}^{\ast }(i)\mathcal{X}_{n}^{\ast }(j)$ and $\Delta
_{n}$ $:=$ $\max_{i,j}|\sigma _{n}^{2}(i,j)$ $-$ $\sigma _{n}^{\ast 2}(i,j)|$%
. By Theorem 2 in \cite{Chernozhukov_etal2015},%
\begin{equation}
\delta _{n}^{\ast }\lesssim \Delta _{n}^{1/3}\times \left\{ 1\vee \ln
(p/\Delta _{n})\right\} ^{2/3}\leq \Delta _{n}^{1/3}\times \left\{ 1+\ln
(p)+\left\vert \ln \Delta _{n}\right\vert \right\} ^{2/3}:=d(n,p).
\label{dn*}
\end{equation}%
It suffices to prove $\Delta _{n}\lesssim \chi _{n}/b_{n}$. Subsequently (%
\ref{rhotil}) with (\ref{dn*}) yields (\ref{rex1_til}).

We now prove $\Delta _{n}\lesssim \chi _{n}/b_{n}$. Using %
\citet[Corollary]{Davydov1968}'s bound under geometric mixing, mutual
independence and $\mathbb{E}\varepsilon _{l}^{2}$ $=$ $1$, and Lyapunov's
inequality, we deduce for some $q$ $>$ $2$ and $\omega $ $\in $ $(0,1)$,%
\footnote{%
A sharp covariance bound, modulo constants, due to \cite{Rio1993} will not
improve the rate at which $\Delta _{n}$ $\rightarrow $ $0$.} 
\begin{eqnarray*}
\Delta _{n} &=&\max_{i,j}\left\vert \frac{1}{n}\sum_{l=1}^{\mathcal{N}%
_{n}}\sum_{s=(l-1)b_{n}+1}^{lb_{n}}\sum_{t\notin (l-1)b_{n}+1}^{lb_{n}}%
\mathbb{E}x_{i,s}x_{j,t}\right\vert \\
&\leq &\max_{i,t}\left\Vert x_{i,t}\right\Vert _{q}\times \max_{i}\left\vert 
\frac{1}{\mathcal{N}_{n}}\sum_{l=1}^{\mathcal{N}_{n}}\frac{1}{b_{n}}%
\sum_{s=(l-1)b_{n}+1}^{lb_{n}}\sum_{t\notin (l-1)b_{n}+1}^{lb_{n}}\alpha
_{i,j,n}(\left\vert s-t\right\vert )^{1-2/q}\right\vert \\
&\leq &\max_{i,t}\left\Vert x_{i,t}\right\Vert _{q}\times \max_{1\leq l\leq 
\mathcal{N}_{n}}\left\vert \frac{1}{b_{n}}\sum_{s=(l-1)b_{n}+1}^{lb_{n}}%
\sum_{t\notin (l-1)b_{n}+1}^{lb_{n}}\omega ^{\left\vert s-t\right\vert
(1-2/q)}\right\vert .
\end{eqnarray*}%
It is straightforward to verify $\sum_{s=(l-1)b_{n}+1}^{lb_{n}}\sum_{t\notin
(l-1)b_{n}+1}^{lb_{n}}\omega ^{\left\vert s-t\right\vert (1-2/q)}$ $=$ $O(1)$
for each $l$ $\in $ $\{1,...,\mathcal{N}_{n}\}$ since the summand index sets 
$s$ $\in $ $\mathfrak{B}_{l}$ and $t$ $\in $ $\mathfrak{B}_{l}^{c}$\ are
mutually exclusive. This proves $\Delta _{n}\lesssim
\max_{i,t}||x_{i,t}||_{q}/b_{n}$. Thus the proof is complete since $%
\max_{i,t}||x_{i,t}||_{q}$ $\lesssim $ $q\chi _{n}$ from (\ref{Exq}). $%
\mathcal{QED}$.\medskip \newline
\noindent \textbf{Proof of Lemma \ref{lm:ex_phys}}. \medskip \newline
\textbf{Claim (}$a$\textbf{)}. Assumptions 2.1-2.3 in \cite{ZhangCheng2018}
\lbrack ZC] hold by Assumption \ref{assume:ex_phys}. Hence (\ref{rhonM})
holds by Theorem 2.1 in ZC when $q$ $\geq $ $2$.\medskip \newline
\textbf{Claim (}$b$\textbf{)}. Recall $y_{i,t}^{(\mathcal{M})}$ $:=$ $%
x_{i,t}^{(\mathcal{M})}$ $-$ $\mathbb{E}x_{i,t}^{(\mathcal{M})}$ with $%
x_{i,t}^{(\mathcal{M})}$ $:=$ $x_{i,t}\mathcal{I}_{\left\vert
x_{i,t}\right\vert \leq \mathcal{M}_{n}}$ $+$ $\mathcal{M}_{n}\mathcal{I}%
_{\left\vert x_{i,t}\right\vert >\mathcal{M}_{n}}$, and $\mathcal{S}_{n,l}^{(%
\mathcal{M})}(i):=\sum_{t=(l-1)b_{n}+1}^{lb_{n}}y_{i,t}^{(\mathcal{M})}$.
Also $\mathcal{X}_{n}^{(\mathcal{M})}(i)$ $=$ $1/\sqrt{n}%
\sum_{t=1}^{n}y_{i,t}^{(\mathcal{M})}$ with blocked version $\mathcal{X}%
_{n}^{(\mathcal{M})\ast }(i)$ $=$ $1/\sqrt{\mathcal{N}_{n}}\sum_{l=1}^{%
\mathcal{N}_{n}}\varepsilon _{l}\mathcal{S}_{n,l}^{(\mathcal{M})}(i)/\sqrt{%
b_{n}}$.

Let $\{\boldsymbol{X}_{n}^{(\mathcal{M})}(i),\boldsymbol{X}_{n}^{(\mathcal{M}%
)\ast }(i)\}$ be Gaussian processes with marginals $\boldsymbol{X}_{n}^{(%
\mathcal{M})}(i)$ $\sim $ $N(0,\mathbb{E}\mathcal{X}_{n}^{(\mathcal{M}%
)}(i)^{2})$ and $\boldsymbol{X}_{n}^{(\mathcal{M})\ast }(i)$ $\sim $ $N(0,%
\mathbb{E}\mathcal{X}_{n}^{(\mathcal{M})\ast }(i)^{2})$. We will bound for
some $\{c_{n}(p),d_{n}(p)\}$ the following Gaussian approximation and
Gaussian-to-Gaussian comparison, 
\begin{eqnarray}
&&\rho _{n}^{(\mathcal{M})\ast \ast }:=\sup_{z\geq 0}\left\vert \mathbb{P}%
\left( \max_{i}\left\vert \mathcal{X}_{n}^{(\mathcal{M})\ast }(i)\right\vert
\leq z\right) -\mathbb{P}\left( \max_{i}\left\vert \boldsymbol{X}_{n}^{(%
\mathcal{M})\ast }(i)\right\vert \leq z\right) \right\vert \lesssim c_{n}(p)%
\text{ \ \ \ }  \label{rho*M_til} \\
&&\delta _{n}^{(\mathcal{M})\ast }:=\sup_{z\geq 0}\left\vert \mathbb{P}%
\left( \max_{i}\left\vert \boldsymbol{X}_{n}^{(\mathcal{M})}(i)\right\vert
\leq z\right) -\mathbb{P}\left( \max_{i}\left\vert \boldsymbol{X}_{n}^{(%
\mathcal{M})\ast }(i)\right\vert \leq z\right) \right\vert \lesssim d_{n}(p),
\label{d*M}
\end{eqnarray}%
thus $\rho _{n}^{\ast (\mathcal{M})}$ $\lesssim $ $c_{n}(p)$ $\vee $ $%
d_{n}(p)$.\medskip \newline
\textbf{Step 1: (\ref{rho*M_til}).} It suffices to will prove $\mathcal{\{}%
\varepsilon _{l}\mathcal{S}_{n,l}^{(\mathcal{M})}(i)/\sqrt{b_{n}}\}_{l=1}^{%
\mathcal{N}_{n}}$ satisfies Assumptions 2.1-2.3 in ZC. Then by their Theorem
2.1%
\begin{eqnarray*}
\rho _{n}^{(\mathcal{M})\ast \ast } &\lesssim &\frac{\mathcal{\tilde{E}}%
_{n}^{1/2}\tilde{l}_{n}^{7/8}}{n^{1/8}}+\tilde{\gamma}_{n}+\left( \frac{%
n^{1/8}}{\mathcal{\tilde{E}}_{n}^{1/2}\tilde{l}_{n}^{3/8}}\right)
^{q/(q+1)}\left( \sum_{i=1}^{p}\check{\Theta}_{i,n}^{(\mathcal{M}%
)(q/2)}(m)\right) ^{1/(q+1)} \\
&&+\check{\Xi}_{n}^{(\mathcal{M})1/3}(m)\times \left( 1\vee \ln \left( p/%
\check{\Xi}_{n}(m)\right) \right) ^{2/3}:=c_{n}(p),
\end{eqnarray*}%
where $\tilde{l}_{n}$ $:=$ $l_{n}(p,\gamma )$ $=$ $\ln (pn/\tilde{\gamma}%
_{n})\vee 1$; $\mathcal{\tilde{E}}_{n}$ $\geq $ $1$ and $\tilde{\gamma}_{n}$ 
$\in $ $(0,1)$ satisfy%
\begin{equation*}
\frac{n^{3/8}}{\mathcal{\tilde{E}}_{n}^{1/2}\tilde{l}_{n}^{5/8}}\geq
K\left\{ \mathcal{M}_{n}h^{-1}\left( n/\tilde{\gamma}_{n}\right) \vee \tilde{%
l}_{n}^{1/2}\right\} ;
\end{equation*}%
and $\check{\Theta}_{i,n}^{(\mathcal{M})(q)}(m)$ $:=$ $\sum_{l=m}^{\infty }%
\check{\Theta}_{i,t}^{(\mathcal{M})(q)}(l)$ and $\check{\Xi}_{n}^{(\mathcal{M%
})}(m)$ $:=$ $\max_{i}\sum_{l=m}^{\infty }l\tilde{\theta}_{i,n}^{(\mathcal{M}%
)(2)}(l)$ with 
\begin{equation*}
\tilde{\theta}_{i,n}^{(\mathcal{M})(q)}(m):=\max_{1\leq l\leq \mathcal{N}%
_{n}}\left\Vert \varepsilon _{l}\frac{1}{\sqrt{b_{n}}}%
\sum_{t=(l-1)b_{n}+1}^{lb_{n}}x_{i,t}^{(\mathcal{M})}-\varepsilon
_{l}^{\prime }(m)\frac{1}{\sqrt{b_{n}}}%
\sum_{t=(l-1)b_{n}+1}^{lb_{n}}x_{i,t}^{(\mathcal{M})\prime }(m)\right\Vert
_{q}.
\end{equation*}%
We show below that $\mathcal{\tilde{E}}_{n}$ $=$ $\mathcal{E}_{n}$, $\tilde{%
\gamma}_{n}$ $=$ $\gamma _{n}$ and therefore $\tilde{l}_{n}$ $=$ $l_{n}$
suffice.\medskip

\textbf{Assumption 2.1 in ZC}: We need for some sequence $\{\tilde{\chi}%
_{n}\},$ $\tilde{\chi}_{n}$ $\geq $ $1$, 
\begin{equation}
(i)\text{ }\max_{i,l}\left\Vert \frac{\varepsilon _{l}\mathcal{S}_{n,l}^{(%
\mathcal{M})}(i)}{\sqrt{b_{n}}}\right\Vert _{4}\leq K\text{ \ and \ }(ii)%
\text{ }\max_{l}\mathbb{E}h\left( \frac{\left\vert \varepsilon _{l}\mathcal{S%
}_{n,l}^{(\mathcal{M})}(i)\right\vert }{\tilde{\chi}_{n}\sqrt{b_{n}}}\right)
\leq K.  \label{A21}
\end{equation}%
The first bound (\ref{A21}).$(i)$ holds by mutual independence of $%
(\varepsilon _{l},\mathcal{S}_{n,l}^{(\mathcal{M})}(i))$, $|\varepsilon
_{l}| $ $\leq $ $c$ $a.s.$, physical dependence Assumption \ref%
{assume:ex_phys}.c, and thus Theorems 1 and 2($i$) in \cite{Wu2005}: 
\begin{equation*}
\max_{i,l}\left\Vert \frac{\varepsilon _{l}\mathcal{S}_{n,l}^{(\mathcal{M}%
)}(i)}{\sqrt{b_{n}}}\right\Vert _{4}\leq c\max_{i,l}\left\Vert \frac{1}{%
\sqrt{b_{n}}}\sum_{t=(l-1)b_{n}+1}^{lb_{n}}y_{i,t}^{(\mathcal{M}%
)}\right\Vert _{4}\lesssim \sum_{m=0}^{\infty }\max_{i}\theta _{i,n}^{%
\mathcal{M}(4)}(m)=O(1).
\end{equation*}

Next, set $\tilde{\chi}_{n}$ $=$ $4\max \{\chi _{n},\mathcal{M}_{n}\}\sqrt{%
b_{n}}$. By $h$ convex and strictly increasing, mutual independence, $%
|\varepsilon _{l}|$ $\leq $ $c$ $a.s.$, and Jensen's inequality, 
\begin{eqnarray*}
\max_{t}\mathbb{E}h\left( \frac{\left\vert \varepsilon _{l}\mathcal{S}%
_{n,l}^{(\mathcal{M})}(i)\right\vert }{\tilde{\chi}_{n}\sqrt{b_{n}}}\right)
&\leq &\max_{l}\frac{1}{b_{n}}\sum_{t=(l-1)b_{n}+1}^{lb_{n}}\mathbb{E}%
h\left( \frac{\left\vert \varepsilon _{l}y_{i,t}^{(\mathcal{M})}\right\vert 
}{4\max \left\{ \chi _{n},\mathcal{M}_{n}\right\} }\right) \\
&\leq &\frac{1}{2}\max_{l}\frac{1}{b_{n}}\sum_{t=(l-1)b_{n}+1}^{lb_{n}}%
\mathbb{E}h\left( \frac{c\left\vert x_{i,t}\right\vert }{\chi _{n}}\right) \\
&&\text{ \ \ \ }+\frac{1}{2}\max_{l}\frac{1}{b_{n}}%
\sum_{t=(l-1)b_{n}+1}^{lb_{n}}\mathbb{E}h\left( \frac{c\left\vert \mathcal{M}%
_{n}\mathcal{I}_{\left\vert x_{i,t}\right\vert >\mathcal{M}_{n}}\right\vert 
}{\mathcal{M}_{n}}\right) \\
&& \\
&\leq &\frac{1}{2}\max_{t}\mathbb{E}h\left( \frac{c\left\vert
x_{i,t}\right\vert }{\chi _{n}}\right) +\frac{h\left( c\right) }{2}\leq K.
\end{eqnarray*}%
The final inequality follows from Assumption \ref{assume:ex_phys}.a. Hence $%
\max_{t}\mathbb{E}h(|\varepsilon _{l}\mathcal{S}_{n,l}^{(\mathcal{M})}(i)|/[%
\tilde{\chi}_{n}\sqrt{b_{n}}])$ $\leq $ $K$ which is (\ref{A21}).$(ii)$%
.\medskip

\textbf{Assumption 2.2 in ZC}: We need to show there exist $\mathcal{\tilde{E%
}}_{n}$ $>$ $0$ and $\tilde{\gamma}_{n}$ $\in $ $(0,1)$ such that $n^{3/8}/[%
\mathcal{\tilde{E}}_{n}^{1/2}\tilde{l}_{n}^{5/8}]$ $\geq $ $K\{\tilde{\gamma}%
_{n}h^{-1}\left( n/\tilde{\gamma}_{n}\right) \vee \tilde{l}_{n}^{1/2}\}$.
Use Assumption \ref{assume:ex_phys}.b with $\mathcal{\tilde{E}}_{n}$ $=$ $%
\mathcal{E}_{n}$ and $\tilde{\gamma}_{n}$ $=$ $\gamma _{n}$, and therefore $%
\tilde{l}_{n}$ $=$ $l_{n}$.\medskip

\textbf{Assumption 2.3 in ZC}: Define%
\begin{equation*}
\sigma _{n}^{(\mathcal{M})\ast 2}(i):=\mathbb{E}\mathcal{X}_{n}^{(\mathcal{M}%
)\ast }(i)^{2}=\mathbb{E}\left( \frac{1}{\sqrt{\mathcal{N}_{n}}}\sum_{l=1}^{%
\mathcal{N}_{n}}\frac{\varepsilon _{l}\mathcal{S}_{n,l}^{(\mathcal{M})}(i)}{%
\sqrt{b_{n}}}\right) ^{2}.
\end{equation*}%
We need to show 
\begin{equation}
0<\min_{i}\sigma _{n}^{(\mathcal{M})\ast 2}(i)\leq \max_{i}\sigma _{n}^{(%
\mathcal{M})\ast 2}(i)\leq K\text{ for each }n  \label{A23_1}
\end{equation}%
and that $\{\varepsilon _{l}\mathcal{S}_{n,l}^{(\mathcal{M})}(i)/\sqrt{b_{n}}%
\}$ are uniformly $\mathcal{L}_{3}$-physical dependent with coefficients $%
\{\theta _{i,n}^{\mathcal{M\ast }(3)}(m)\}$, 
\begin{equation}
\max_{i}\sum_{l=1}^{\infty }l\theta _{i,n}^{\mathcal{M\ast }(3)}(m)\leq K.
\label{A23_2}
\end{equation}%
Consider (\ref{A23_1}). Assumption \ref{assume:ex_phys}.c gives the lower
bound. The upper bound holds by (\ref{A23_2}), and Theorems 1 and 2($i$) in 
\cite{Wu2005}.

It remains to verify (\ref{A23_2}). By mutual independence of $(\varepsilon
_{l},\mathcal{S}_{n,l}^{(\mathcal{M})}(i))$,%
\begin{eqnarray}
\theta _{i,n}^{\mathcal{M\ast }(3)}(m) &\leq &2\max_{l}\left\Vert \frac{1}{%
\sqrt{b_{n}}}\sum_{t=(l-1)b_{n}+1}^{lb_{n}}\varepsilon _{l}\left( y_{i,t}^{(%
\mathcal{M})}-y_{i,t}^{(\mathcal{M})\prime }(m)\right) \right\Vert _{3} 
\notag \\
&=&2\max_{l}\left\{ \mathbb{E}\left[ \mathbb{E}\left( \left\vert \frac{1}{%
\sqrt{b_{n}}}\sum_{t=(l-1)b_{n}+1}^{lb_{n}}\varepsilon _{l}\left( y_{i,t}^{(%
\mathcal{M})}-y_{i,t}^{(\mathcal{M})\prime }(m)\right) \right\vert ^{3}|%
\mathcal{Y}_{i,n}\right) \right] \right\} ^{1/3}  \notag \\
&\lesssim &\max_{l}\left\Vert y_{i,t}^{(\mathcal{M})}-y_{i,t}^{(\mathcal{M}%
)\prime }(m)\right\Vert _{3}\lesssim \theta _{i,n}^{\mathcal{M}(3)}(m),
\label{th_th}
\end{eqnarray}%
because by independence of $\varepsilon _{l}$, $|\varepsilon _{l}|$ $\leq $ $%
c$ $a.s.$, and Theorems 1 and 2($i$) in \cite{Wu2005}, 
\begin{eqnarray*}
&&\max_{l}\left\Vert \frac{1}{\sqrt{b_{n}}}\sum_{t=(l-1)b_{n}+1}^{lb_{n}}%
\varepsilon _{l}\left( y_{i,t}^{(\mathcal{M})}-y_{i,t}^{(\mathcal{M})\prime
}(m)\right) |\mathcal{Y}_{i,n}\right\Vert _{3} \\
&&\text{ \ \ \ \ \ \ \ }\leq \sum_{r=0}^{\infty }\max_{l}\left\Vert
\varepsilon _{l}\left( y_{i,t}^{(\mathcal{M})}-y_{i,t}^{(\mathcal{M})\prime
}(m)\right) -\varepsilon _{l}^{\prime }(r)\left( y_{i,t}^{(\mathcal{M}%
)}-y_{i,t}^{(\mathcal{M})\prime }(m)\right) \right\Vert _{3} \\
&&\text{ \ \ \ \ \ \ \ }=\max_{l}\left\Vert \left( \varepsilon
_{l}-\varepsilon _{l+1}\right) \left( y_{i,t}^{(\mathcal{M})}-y_{i,t}^{(%
\mathcal{M})\prime }(m)\right) \right\Vert _{3}\leq 2c\max_{l}\left\Vert
y_{i,t}^{(\mathcal{M})}-y_{i,t}^{(\mathcal{M})\prime }(m)\right\Vert _{3}.
\end{eqnarray*}%
Now exploit Assumption \ref{assume:ex_phys}.c with (\ref{th_th}) to yield (%
\ref{A23_2}).\medskip \newline
\textbf{Step 2: (\ref{d*M}). }Define $\sigma _{n}^{(\mathcal{M})2}(i,j)$ $:=$
$\mathbb{E}\mathcal{X}_{n}^{(\mathcal{M})}(i)\mathcal{X}_{n}^{(\mathcal{M}%
)}(j)$, $\sigma _{n}^{(\mathcal{M})\ast 2}(i,j)$ $:=$ $\mathbb{E}\mathcal{X}%
_{n}^{(\mathcal{M})\ast }(i)\mathcal{X}_{n}^{(\mathcal{M})\ast }(j)$, and $%
\Delta _{n}^{(\mathcal{M})}$ $:=$ $\max_{i,j}|\sigma _{n}^{(\mathcal{M}%
)2}(i,j)$ $-$ $\sigma _{n}^{(\mathcal{M})\ast 2}(i,j)|$. We will prove $%
\Delta _{n}^{(\mathcal{M})}$ $=$ $O(\xi _{n}/b_{n})$, hence $\delta _{n}^{(%
\mathcal{M})\ast }$ $\lesssim $ $\left( \xi _{n}/b_{n}\right) ^{1/3}\ln
(p)^{2/3}$ by Lemma C.5 in \cite{Chen_Kato_2019} 
\citep[cf.][Theorem
2]{Chernozhukov_etal2015}.

By construction, mutual independence and $\mathbb{E}\varepsilon _{l}^{2}$ $=$
$1$,%
\begin{eqnarray}
\Delta _{n}^{(\mathcal{M})} &=&\max_{i,j}\left\vert \sigma _{n}^{(\mathcal{M}%
)2}(i,j)-\sigma _{n}^{(\mathcal{M})\ast 2}(i,j)\right\vert  \notag \\
&=&\max_{i,j}\left\vert \frac{1}{n}\sum_{l=1}^{\mathcal{N}%
_{n}}\sum_{s=(l-1)b_{n}+1}^{lb_{n}}\sum_{t\notin (l-1)b_{n}+1}^{lb_{n}}%
\mathbb{E}y_{i,s}^{(\mathcal{M})}y_{j,t}^{(\mathcal{M})}\right\vert .\text{
\ \ \ \ \ \ \ \ \ \ \ \ \ \ \ }  \label{D=}
\end{eqnarray}%
Now define $\mathcal{F}_{t}$ $:=$ $\sigma (\epsilon _{t},\epsilon
_{t-1},\ldots )$ and a projection operator $\mathbb{P}_{t}y_{i,t}^{(\mathcal{%
M})}$ $:=$ $\mathbb{E}_{\mathcal{F}_{t}}y_{i,t}^{(\mathcal{M})}$ $-$ $%
\mathbb{E}_{\mathcal{F}_{t-1}}y_{i,t}^{(\mathcal{M})}$. Then $y_{i,t}^{(%
\mathcal{M})}$ $=$ $\sum_{l=0}^{\infty }\mathbb{P}_{t-l}y_{i,t}^{(\mathcal{M}%
)}$. Therefore by the martingale difference property of $\mathbb{P}%
_{t-l}y_{i,t}^{(\mathcal{M})}$, and triangle and Cauchy-Schwartz
inequalities, if $s$ $\leq $ $t$%
\begin{eqnarray*}
\left\vert \mathbb{E}y_{i,s}^{(\mathcal{M})}y_{j,t}^{(\mathcal{M}%
)}\right\vert &\leq &\left\vert \mathbb{E}\sum_{m_{1}=0}^{\infty }\mathbb{P}%
_{s-m_{1}}y_{i,s}^{(\mathcal{M})}\sum_{m_{2}=0}^{\infty }\mathbb{P}%
_{t-m_{2}}y_{j,t}^{(\mathcal{M})}\right\vert \\
&\leq &\sum_{m=0}^{\infty }\left\vert \mathbb{E}\left( \mathbb{P}%
_{s-m}y_{i,s}^{(\mathcal{M})}\right) \left( \mathbb{P}_{s-m}y_{j,t}^{(%
\mathcal{M})}\right) \right\vert \leq \sum_{m=0}^{\infty }\left\Vert \mathbb{%
P}_{s-m}y_{i,s}^{(\mathcal{M})}\right\Vert _{2}\left\Vert \mathbb{P}%
_{s-m}y_{j,t}^{(\mathcal{M})}\right\Vert _{2},
\end{eqnarray*}%
and if $s$ $>$ $t$ then $|\mathbb{E}y_{i,s}^{(\mathcal{M})}y_{j,t}^{(%
\mathcal{M})}|$ $\leq $ $\sum_{m=0}^{\infty }||\mathbb{P}_{t-m}y_{i,s}^{(%
\mathcal{M})}||_{2}$ $\times $ $||\mathbb{P}_{t-m}y_{j,t}^{(\mathcal{M}%
)}||_{2}$. Use Theorem 1 in \cite{Wu2005} to yield $||\mathbb{P}%
_{t-m}y_{i,t}^{(\mathcal{M})}||_{2}$ $\leq $ $\theta _{i,n}^{\mathcal{M}%
(2)}(m)$. Thus%
\begin{equation}
\left\vert \mathbb{E}y_{i,s}^{(\mathcal{M})}y_{j,t}^{(\mathcal{M}%
)}\right\vert \leq \sum_{m=0}^{\infty }\theta _{i,n}^{\mathcal{M}%
(2)}(m)\times \theta _{j,n}^{\mathcal{M}(2)}(m+\left\vert t-s\right\vert ).
\label{|Eyy|}
\end{equation}%
Combine (\ref{D=}) with (\ref{|Eyy|}) to deduce%
\begin{eqnarray*}
\Delta _{n}^{(\mathcal{M})} &\leq &\max_{i,j}\left\vert \frac{1}{n}%
\sum_{l=1}^{\mathcal{N}_{n}}\sum_{s=(l-1)b_{n}+1}^{lb_{n}}\sum_{t\notin
(l-1)b_{n}+1}^{lb_{n}}\left\{ \sum_{m=0}^{\infty }\theta _{i,n}^{\mathcal{M}%
(2)}(m)\theta _{j,n}^{\mathcal{M}(2)}(m+\left\vert t-s\right\vert )\right\}
\right\vert \\
&=&\max_{i,j}\left\vert \sum_{m=0}^{\infty }\theta _{i,n}^{\mathcal{M}%
(2)}(m)\left\{ \frac{1}{n}\sum_{l=1}^{\mathcal{N}_{n}}%
\sum_{s=(l-1)b_{n}+1}^{lb_{n}}\sum_{t\notin (l-1)b_{n}+1}^{lb_{n}}\theta
_{j,n}^{\mathcal{M}(2)}(m+\left\vert t-s\right\vert )\right\} \right\vert .
\end{eqnarray*}%
Finally, it is straightforward to verify%
\begin{eqnarray*}
&&\frac{1}{n}\sum_{l=1}^{\mathcal{N}_{n}}\sum_{s=(l-1)b_{n}+1}^{lb_{n}}%
\sum_{t\notin (l-1)b_{n}+1}^{lb_{n}}\theta _{j,n}^{\mathcal{M}%
(2)}(m+\left\vert t-s\right\vert ) \\
&&\text{ \ \ \ \ \ \ \ \ \ \ \ }\leq \frac{1}{\mathcal{N}_{n}}\sum_{l=1}^{%
\mathcal{N}_{n}}\left( \frac{1}{b_{n}}\sum_{s=(l-1)b_{n}+1}^{lb_{n}}%
\sum_{t=1}^{(l-1)b_{n}}\theta _{j,n}^{\mathcal{M}(2)}(m+\left\vert
t-s\right\vert )\right) \\
&&\text{ \ \ \ \ \ \ \ \ \ \ \ \ \ \ \ \ \ \ \ \ \ \ \ }+\frac{1}{\mathcal{N}%
_{n}}\sum_{l=1}^{\mathcal{N}_{n}}\left( \frac{1}{b_{n}}%
\sum_{s=(l-1)b_{n}+1}^{lb_{n}}\sum_{t=lb_{n}+1}^{n}\theta _{j,n}^{\mathcal{M}%
(2)}(m+\left\vert t-s\right\vert )\right) \\
&& \\
&&\text{ \ \ \ \ \ \ \ \ \ \ \ }=\frac{1}{\mathcal{N}_{n}}\sum_{l=1}^{%
\mathcal{N}_{n}}\left( \frac{1}{b_{n}}\sum_{r=1}^{lb_{n}}r\theta _{j,n}^{%
\mathcal{M}(2)}(m+r)\right) \\
&&\text{ \ \ \ \ \ \ \ \ \ \ \ \ \ \ \ \ \ \ \ \ \ \ \ }+\frac{1}{\mathcal{N}%
_{n}}\sum_{l=1}^{\mathcal{N}_{n}}\left( \frac{1}{b_{n}}%
\sum_{r=1}^{n-(l-1)b_{n}-1}r\theta _{j,n}^{\mathcal{M}(2)}(m+r)\right) .
\end{eqnarray*}%
Under Assumption \ref{assume:ex_phys}.c and Lyapunov's inequality $%
\max_{j}\max_{m\geq 0}1/b_{n}\sum_{r=1}^{\infty }r\theta _{j,n}^{\mathcal{M}%
(2)}(m$ $+$ $r)$ $=$ $O\left( \xi _{n}/b_{n}\right) $. Hence%
\begin{equation*}
\max_{j}\max_{m\geq 0}\frac{1}{n}\sum_{l=1}^{\mathcal{N}_{n}}%
\sum_{s=(l-1)b_{n}+1}^{lb_{n}}\sum_{t\notin (l-1)b_{n}+1}^{lb_{n}}\theta
_{j,n}^{\mathcal{M}(2)}(m+\left\vert t-s\right\vert )=O\left( \xi
_{n}/b_{n}\right) ,
\end{equation*}%
and therefore $\Delta _{n}^{(\mathcal{M})}$ $\lesssim $ $b_{n}^{-1}\max_{i}|%
\sum_{m=0}^{\infty }\theta _{i,n}^{\mathcal{M}(2)}(m)|$ $=$ $O\left( \xi
_{n}/b_{n}\right) $ as claimed$.$ $\mathcal{QED}$.

\bibliographystyle{cas-model2-names.bst}
\bibliography{refs_max_ineq}
\clearpage

\end{document}